\documentclass[12pt,a4paper]{amsart}
\usepackage{amssymb,verbatim}

\theoremstyle{plain}
\newtheorem{theorem}{Theorem}
\newtheorem{proposition}[theorem]{Proposition}
\newtheorem{lemma}[theorem]{Lemma}

\theoremstyle{definition}
\newtheorem*{remark}{Remark}
\newtheorem*{definition}{Definition}

\newcommand{\prz}{|X(t,\omega)-X(s,\omega)|}
\newcommand{\il}{\bigcap_{\substack{|t-s|<1/m\\s,t\in[a,b]}}}
\newcommand{\omo}{\Omega_0}

\begin{document}
\title[Stochastic integrals and Doob-Meyer decomposition]
{Quantum stochastic integrals and Doob-Meyer decomposition}
\author[A.\L uczak]{Andrzej \L uczak}
\address{Faculty of Mathematics \\
        \L\'od\'z University \\
        ul. S. Banacha 22 \\
        90-238 \L\'od\'z, Poland}
\email{anluczak@math.uni.lodz.pl}
\thanks{Work supported by KBN grant 2PO3A 03024}
\keywords{Quantum stochastic integrals, Doob-Meyer decomposition,
quantum martingales, von Neumann algebras} \subjclass{Primary:
81S25; Secondary: 46L53, 60G44, 60H05}
\date{}
\begin{abstract}
We show that for a quantum $L^p$-martingale $(X(t))$, $p>2$, there
exists a Doob-Meyer decomposition of the submartingale
$(|X(t)|^2)$. A noncommutative counterpart of a classical process
continuous with probability one is introduced, and a quantum
stochastic integral of such a process with respect to an
\mbox{$L^p$-mar}\-tingale, $p>2$, is constructed. Using this
construction, the uniqueness of the Doob-Meyer decomposition for a
quantum martingale `continuous with probability one' is proved,
and explicit forms of this decomposition and the quadratic
variation process for such a martingale are obtained.
\end{abstract}
\maketitle

\section*{Introduction}
In the existing theories of quantum stochastic integration we face
a problem which is similar to that described in \cite[p. 148]{KS}
as follows: ``We know all to well that it is one thing to develop
a theory of integration in some reasonable generality and a
completely different task to \emph{compute} the integral in any
specific case of interest''. In quantum stochastic integration we
are concerned not so much with computing the integral but rather
with numerous important examples which do not fit into the nice
theory that we have at our disposal. The origin of this problem
lies in rather narrow classes of `theoretically admissible'
integrands. Indeed, if for example $(X(t))$ is an $L^2$-martingale
then the integral $\int f(t)\,dX(t)$ is in general defined for
adapted processes $f$ satisfying pretty strong conditions such as
e.g. being \emph{norm} limits of simple processes. On the other
hand it looks quite reasonable to define a stochastic integral in
some natural way in many concrete situations making it possible to
integrate a broader class of processes. This approach has already
been taken in \cite{BW4,BW5} in several cases, where in particular
it is shown how integration with respect to a quantum random time
can be performed, or how can one integrate predictable processes.
In the first part of the paper we follow the same idea and
construct a stochastic integral of `continuous with probability
one' noncommutative stochastic process with respect to an
$L^p$-martingale for $p>2$. The second part is devoted to the
problem of a Doob-Meyer decomposition of the submartingale
$(|X(t)|^2)$, where $(X(t))$ is an $L^p$-martingale for $p>2$.
This problem has a long history, cf. for example
\cite{BGW,BSW1,BSW2,BW1,BW2,BW3,BW4,BW5}. We show that such a
decomposition always exists, and is unique for a martingale
`continuous with probability one'. In this case we give also
explicit forms of the quadratic variation process of the
martingale and the Doob-Meyer decomposition, using the
construction of the integral given in the first part of the paper.
As seen from the above, the notion of `continuity with probability
one' for a noncommutative stochastic process plays an important
role in our considerations, and we explain how it can be
generalised from the classical context to the noncommutative one.

\section{Preliminaries and notation}
A \emph{noncommutative stochastic base} which is a basic object of
our considerations consists of the following elements: a von
Neumann algebra $\mathcal{A}$ acting on a Hilbert space
$\mathcal{H}$, a normal faithful unital trace $\tau$ on
$\mathcal{A}$, a filtration $(\mathcal{A}_t\colon t \in
[0,+\infty))$, which is an increasing ($s \leqslant t$  implies
$\mathcal{A}_s \subset \mathcal{A}_t$) family of von Neumann
subalgebras of $\mathcal{A}$ such that
$\mathcal{A}=\mathcal{A}_{\infty} = \left( \bigcup_{t \geqslant 0}
\mathcal{A}_t \right)''$ and $\mathcal{A}_s = \bigcap_{t>s}
\mathcal{A}_t$ (right-continuity). Moreover, for each $t$, there
exists a normal conditional expectation $\mathbb{E}_t$ from
$\mathcal{A}$ onto $\mathcal{A}_t$ such that $\tau \circ
\mathbb{E}_t = \tau$.

For each $t\in[0,+\infty]$ we write $L^p(\mathcal{A}_t)$ for the
non-commutative Lebesgue space associated with $\mathcal{A}_t$ and
$\tau$. The theory of such spaces is described e.g. in \cite{Y};
for our purposes we only recall that $L^p(\mathcal{A})$
(respectively $L^p(\mathcal{A}_t)$) consists of densely defined
operators on $\mathcal{H}$, affiliated to $\mathcal{A}$, and that
$L^p(\mathcal{A})$ is completion of $\mathcal{A}$ with respect to
the norm
\[
\|X\|_p = [\tau (|X|^p)]^{1\slash p};
\]
moreover, for $a \in \mathcal{A}$, $X \in L^p(\mathcal{A})$ the
operators $aX$ and $Xa$ belong to $L^p(\mathcal{A})$. For each $t$
the conditional expectation $\mathbb{E}_t$ extends to a projection
of norm one from $L^p(\mathcal{A})$ onto
$L^p(\mathcal{A}_t)$---for which we use the same notation. Notice
that the conditional expectation, being a bounded operator on
$L^p(\mathcal{A})$, is weakly continuous. Since the conditional
expectation is completely positive we have
\[
 \mathbb{E}_t|x|^2=\mathbb{E}_t(x^*x)\geqslant\mathbb{E}_tx^*
 \mathbb{E}_tx=|\mathbb{E}_tx|^2.
\]
The following simple property is often useful. Let $x\in
L^p(\mathcal{A}),\linebreak y\in
L^q(\mathcal{A}),\,p,q\in[0,+\infty],\,\frac{1}{p}+\frac{1}{q}=1$.
Then
\begin{equation}\label{ce}
 \tau((\mathbb{E}_tx)y)=\tau(\mathbb{E}_t((\mathbb{E}_tx)y))
 =\tau((\mathbb{E}_tx)(\mathbb{E}_ty))=\tau(x\mathbb{E}_ty).
\end{equation}

By an $\mathcal{A}$-(respectively $L^p$-) valued \emph{process} we
mean a map from $[0,+\infty)$ into $\mathcal{A}$ (respectively
$L^p(\mathcal{A})$). $\mathcal{A}$-valued processes will be
usually denoted by $f$, $g$ or $(f(t))$ , $(g(t))$, while for
$L^p$-processes we shall use symbols $(X(t))$, $(Y(t))$. An
$\mathcal{A}$-(respectively $L^p$-) valued process $f$
(respectively $X$) is called \emph{adapted} if $f(t) \in
\mathcal{A}_t$ (respectively \linebreak $X(t) \in L^p
(\mathcal{A}_t)$).

An $L^p$-process $(X(t) \colon t \in [0,+\infty))$ is called a
\emph{martingale} if for each $s,t \in [0,\infty)$, $s \leqslant
t$, we have $\mathbb{E}_s X(t) = X(s)$. It follows that a
martingale is an adapted process. If the inequality
$X(s)\leqslant\mathbb{E}_sX(t)$ (resp. $\mathbb{E}_sX(t)\leqslant
X(s)$) holds for $s\leqslant t$, then the process is called a
\emph{submartingale} (resp. \emph{supermartingale}). Let us notice
that according to \cite{B} the martingale $(X(t))$ is
right-continuous in $\|\cdot\|_p$-norm. Moreover, for each
$p\in[0,+\infty]\text{ and }s\leqslant t$ we have
\[
 \|X(s)\|_p=\|\mathbb{E}_sX(t)\|_p\leqslant\|X(t)\|_p.
\]
If $(X(t)\colon t\in[0,+\infty))$ is an $L^p$-martingale, for
$p\geqslant 2$, then \linebreak $(|X(t)|^2\colon t\in[0,+\infty))$
is an $L^{p/2}$-submartingale since for any $s\leqslant t$ the
above-mentioned property of conditional expectation yields
\[
 \mathbb{E}_s|X(t)|^2\geqslant|\mathbb{E}_sX(t)|^2=|X(s)|^2.
\]
The submartingale $(|X(t)|^2\colon t\in[0,+\infty))$ is
right-continuous in \linebreak $\|\cdot\|_{p/2}$-norm. Indeed, we
have
\[
 |X(t)|^2-|X(s)|^2=X(t)^*[X(t)-X(s)]+[X(t)-X(s)]^*X(s),
\]
so using H\"{o}lder's inequality we get
\begin{align*}
 &\||X(t)|^2-|X(s)|^2\|_{p/2}\leqslant
 \|X(t)^*[X(t)-X(s)]\|_{p/2}\\
 +&\|[X(t)-X(s)]^*X(s)\|_{p/2} \leqslant
 \|X(t)\|_p\|X(t)-X(s)\|_p\\+&\|X(t)-X(s)\|_p\|X(s)\|_p,
\end{align*}
which on account of the right-continuity of $(X(t))$ in
$\|\cdot\|_p$-norm shows that for $t\searrow s,\quad |X(t)|^2\to
|X(s)|^2\text{ in }\|\cdot\|_{p/2}$-norm.

Let $(X(t)\colon t\in[0,+\infty))$ be a process, and let
$0\leqslant t_0\leqslant t_1\leqslant\dots\leqslant t_m<+\infty$
be a sequence of points. To simplify the notation we put
\[
 \Delta X(t_k)=X(t_k)-X(t_{k-1}),\qquad k=1,\dots,m.
\]

Let $(X(t)\colon t\in[0,+\infty)),\,(Y(t)\colon t\in[0,+\infty))$
be arbitrary processes, and let $[a,b]$ be a subinterval of
$[0,+\infty)$. For a partition \linebreak
$\theta=\{a=t_0<t_1<\dots<t_m=b\}$ of $[a,b]$ we form left and
right integral sums
\begin{align*}
 &S_{\theta}^l=\sum_{k=1}^m\Delta X(t_k)Y(t_{k-1})\\
 &S_{\theta}^r=\sum_{k=1}^m Y(t_{k-1})\Delta X(t_k).
\end{align*}
If there exist limits (in any sense) of the above sums as $\theta$
refines, we call them respectively the \emph{left} and \emph{right
stochastic integrals} of $(Y(t))$ with respect to $(X(t))$, and
denote
\begin{align*}
 &\lim_{\theta}S_{\theta}^l=\int_a^b dX(t)\,Y(t)\\
 &\lim_{\theta}S_{\theta}^r=\int_a^b Y(t)\,dX(t).
\end{align*}
This notion of integral is a weaker one. Indeed, we could define
the integrals as the limits
\begin{align*}
 &\int_a^b dX(t)\,Y(t)=\lim_{\|\theta\|\to 0}S_{\theta}^l\\
 &\int_a^b Y(t)\,dX(t)=\lim_{\|\theta\|\to 0}S_{\theta}^r,
\end{align*}
where $\|\theta\|$ stands for the mesh of the partition $\theta$.
A definition of this kind is standard in the classical theories of
Riemann-Stieltjes as well as stochastic integral; it is worth
noticing that in noncommutative integration theory, whenever this
Riemann-Stieltjes type integral is considered, its definition
refers to the weaker form of the limit with the refining net of
partitions (cf. \cite{BL,BW4,BW5}). However, in our case we shall
be able to obtain the integral in the stronger sense thus making
it similar to the classical stochastic integral.

Let $(X(k)\colon k=0,1,\dots n)$ be a finite martingale. We have
\begin{equation}\label{e0}
 \begin{aligned}
 &\mathbb{E}_{k-1}|\Delta X(k)|^2=\mathbb{E}_{k-1}(|X(k)|^2-X(k)^*X(k-1)+\\
 -&X(k-1)^*X(k)+|X(k-1)|^2)\\
 =&\mathbb{E}_{k-1}(|X(k)|^2)-(\mathbb{E}_{k-1}X(k)^*)X(k-1)+\\
 -&X(k-1)^*(\mathbb{E}_{k-1}X(k))+|X(k-1)|^2\\
 =&\mathbb{E}_{k-1}|X(k)|^2-|X(k-1)|^2=\mathbb{E}_{k-1}(|X(k)|^2-|X(k-1)|^2),
 \end{aligned}
\end{equation}
by martingale property. From the above we obtain on account of the
$\mathbb{E}_k$-invariance of $\tau$
\begin{equation}\label{e1}
 \begin{aligned}
  &\sum_{k=1}^n\||\Delta X(k)|^2\|_1=\sum_{k=1}^n\tau(|\Delta
  X(k)|^2)\\
  =&\sum_{k=1}^n\tau(\mathbb{E}_{k-1}(|X(k)|^2)
  =\sum_{k=1}^n\tau(\mathbb{E}_{k-1}(|X(k)|^2-|X(k-1)|))\\
  =&\sum_{k=1}^n\tau(|X(k)|^2-|X(k-1)|^2)=\tau(|X(n)|^2)-\tau(|X(0)|^2).
 \end{aligned}
\end{equation}
The equality above gives the obvious estimation
\begin{align*}
 \Big\|&\Big(\sum_{k=1}^n|\Delta X(k)|^2\Big)^{1/2}\Big\|_2=
 \Big(\sum_{k=1}^n\tau|\Delta X(k)|^2\Big)^{1/2}\\
  =&\Big(\sum_{k=1}^n\||\Delta X(k)|^2\|_1\Big)^{1/2}\leqslant\|X(n)\|_2.
\end{align*}
A fundamental result from \cite{PX}---Theorem 2.1---says that the
estimation of this type is valid for each $p>1$. We shall use this
for $p>2$, in which case it has the form: there exists a constant
$\alpha_p$ depending only on $p$, such that for each
$L^p$-martingale $(X(k)\colon k=0,1,\dots,n)$ we have
\begin{equation}\label{PX}
 \Big\|\Big(\sum_{k=1}^n|\Delta
 X(k)|^2\Big)^{1/2}\Big\|_p\leqslant\alpha_p\|X(n)\|_p.
\end{equation}

\section{Continuity of a noncommutative stochastic process}
Let $(X(t,\cdot)\colon t\in[a,b])$ be a stochastic process over a
probability space $(\Omega,\mathcal{F},P)$. Consider the following
condition: for each $\varepsilon>0$ there is
$\Omega_{\varepsilon}\in\mathcal{F}$ with
$P(\Omega_{\varepsilon})>1-\varepsilon$, such that the
trajectories \linebreak
$\{X(\cdot,\omega)\colon\omega\in\Omega_{\varepsilon}\}$ are
equally uniformly continuous. This can be rewritten as:
\[
 \begin{aligned}\label{c}
 &\text{for each }\varepsilon>0\text{ there is }
 \Omega_{\varepsilon}\in\mathcal{F}\text{ with }
 P(\Omega_{\varepsilon})>1-\varepsilon,\\&\text{having the property:}\\
 &\text{for each }\eta>0 \text{ there is } \delta>0 \text{ such
 that for any }\omega\in\Omega_{\varepsilon}\\ &\text{and any }s,t\in[a,b]
 \text{ with }|t-s|<\delta,\\&\text{we have }\prz\leqslant\eta.
 \end{aligned}\tag{*}
\]
If the above condition is satisfied, then the trajectories of the
process are uniformly continuous with probability one. Indeed,
take $\varepsilon=1/n$, and let
$\Omega_{\varepsilon}=\Omega_{1/n}$ be as above. Put
\[
 \omo=\bigcup_{n=1}^{\infty}\Omega_{1/n}.
\]
Then $P(\omo)=1$, and for each $\omega\in\omo\text{ we have
}\omega\in\Omega_{1/n}$ for some $n$, which means that the
trajectory $X(\cdot,\omega)$ is uniformly continuous.

Now let us assume that the trajectories are uniformly continuous
with probability one, and let $\omo=\{\omega\colon X(\cdot,\omega)
\text{ is uniformly continuous}\}$. We have $P(\omo)=1$, and
\begin{align*}
 \omo&=\bigcap_{r=1}^{\infty}\bigcup_{m=1}^{\infty}\il\left\{\omega\colon\prz
 \leqslant\frac{1}{r}\right\}\\&=\bigcap_{r=1}^{\infty}\bigcup_{m=1}^{\infty}\il
 \left\{\omega\in\omo\colon\prz\leqslant\frac{1}{r}\right\}.
\end{align*}
The continuity of the trajectories for $\omega\in\omo$ implies
that
\begin{align*}
 &\il\left\{\omega\in\omo\colon\prz\leqslant\frac{1}{r}\right\}\\=&\bigcap_
 {\substack{|t-s|<1/m\\s,t\in[a,b]\cap\mathbb{Q}}}
 \left\{\omega\in\omo\colon\prz\leqslant\frac{1}{r}\right\},
\end{align*}
where $\mathbb{Q}$ stands for the rational numbers. It follows
that the set
\[
 \bigcap_{\substack{|t-s|<1/m\\s,t\in[a,b]}}
 \left\{\omega\in\omo\colon\prz\leqslant\frac{1}{r}\right\}
\]
is measurable, and for each positive integer $r$ we have
\begin{align*}
 1&=P\Big(\bigcup_{m=1}^{\infty}\il\left\{\omega\colon\prz
 \leqslant\frac{1}{r}\right\}\Big)\\&=\lim_{m\to\infty}P\Big(\bigcap_
 {\substack{|t-s|<1/m\\s,t\in[a,b]}}
 \left\{\omega\in\omo\colon\prz\leqslant\frac{1}{r}\right\}\Big).
\end{align*}
For any $\varepsilon>0$ and positive integer $r$ choose $m_r$ such
that
\[
 P\Big(\bigcap_
 {\substack{|t-s|<1/m_r\\s,t\in[a,b]}}
 \left\{\omega\in\omo\colon\prz\leqslant\frac{1}{r}\right\}\Big)
 >1-\frac{\varepsilon}{2^r},
\]
and put
\[
 \Omega_{\varepsilon}=\bigcap_{r=1}^{\infty}\bigcap_
 {\substack{|t-s|<1/m_r\\s,t\in[a,b]}}
 \left\{\omega\in\omo\colon\prz\leqslant\frac{1}{r}\right\}.
\]
Then $P(\Omega_{\varepsilon})>1-\varepsilon$. For arbitrary fixed
$\eta>0$ let $r_0$ be such that $1/r_0\leqslant\eta$. Put
$\delta=1/m_{r_0}$. For each $\omega_0\in\Omega_{\varepsilon}$ we
have, in particular, that
$\omega_0\in\left\{\omega\in\omo\colon\prz\leqslant\frac{1}{r_0}\right\}$
for any $s,t\in[a,b]$ with \newline$|t-s|<1/m_{r_0}=\delta$, which
means that
\[
 |X(t,\omega_0)-X(s,\omega_0)|\leqslant\frac{1}{r_0}\leqslant\eta,
\]
showing that condition \eqref{c} holds.
\par
We have thus shown the equivalence of uniform continuity of
trajectories of the process with probability one and condition
\eqref{c}. Since in our case the uniform continuity of
trajectories is equivalent to ordinary continuity, condition
\eqref{c} can be treated simply as another definition of the
classical notion of a continuous stochastic process.

Let us observe that condition \eqref{c} can be given the following
form. Denote by $\chi_E$ the indicator function of the set $E$.
Then condition \eqref{c} becomes:\par for each $\varepsilon>0$
there is $\Omega_{\varepsilon}\in\mathcal{F}$ with
$P(\Omega_{\varepsilon})>1-\varepsilon$, having the property: for
each $\eta>0$ there is $\delta>0$ such that for any $s,t\in[a,b]$
with $|t-s|<\delta$, we have
\begin{equation*}
 \begin{aligned}
 \sup_{\omega\in\Omega_{\varepsilon}}\prz&=\sup_{\omega\in\Omega}
 [\prz\chi_{\Omega_{\varepsilon}}(\omega)]\\&=\|[X(t,\cdot)-X(s,\cdot)]
 \chi_{\Omega_{\varepsilon}}\|_{\infty}\leqslant\eta.
 \end{aligned}
\end{equation*}
The above form is essentially \emph{algebraic}, referring only to
the algebra $L^{\infty}(\Omega)$, which becomes clear if we
replace the inequality $P(\Omega_{\varepsilon})>1-\varepsilon$ by
the equivalent inequality
$\int_{\Omega}\chi_{\Omega_{\varepsilon}}\,dP>1-\varepsilon$. Thus
for a noncommutative process $(X(t)\colon t\in[a,b])$ it can be
given either of the following two forms: `right' and `left',
denoted respectively by (R) and (L). \par For each $\varepsilon>0$
there is a projection $e$ in $\mathcal{A}\text{ with
}\tau(e)>1-\varepsilon$, having the property: for each $\eta>0$
there is $\delta>0$ such that for any $s,t\in[a,b]\text{ with
}|t-s|<\delta$, we have
\begin{equation}
 [X(t)-X(s)]e\in\mathcal{A}\text
 {\qquad and\qquad }\|[X(t)-X(s)]e\|_{\infty}\leqslant\eta,\tag{R}
\end{equation}
or
\begin{equation}
 e[X(t)-X(s)]\in\mathcal{A}\text
 {\qquad and\qquad }\|e[X(t)-X(s)]\|_{\infty}\leqslant\eta,\tag{L}
\end{equation}
where $\|\cdot\|_{\infty}$ denotes the norm in the algebra
$\mathcal{A}$. This form of `noncommutative continuity of
trajectories with probability one', in its right version, has
already been considered before in \cite{GL1,GL2}, where it was
given the name of `Segal's uniform continuity', and some theorems
on this continuity were obtained. However, it is easily seen that
the `right Segal's uniform continuity' which appears in the
conclusions of those theorems can be changed to the `left Segal's
uniform continuity', so the results in \cite{GL1,GL2} give in fact
both forms of this continuity. We shall call a process
\emph{uniformly continuous in Segal's sense} if it satisfies both
(R) and (L) conditions. It is obvious that for a selfadjoint
process conditions (R) and (L) are equivalent. Let now
$\varepsilon,\,e\text{ and }\delta$ be as above. For arbitrary
$s,t\in[a,b],\,s<t$ choose points $s=t_0<t_1<\dots<t_m=t$ such
that $\max_{1\leqslant k\leqslant m}(t_k-t_{k-1})<\delta$. Then
\[
 [X(t)-X(s)]e=[X(t)-X(t_{m-1})]e+\cdots+[X(t_1)-X(s)]e,
\]
and since all the summands on the right hand side belong to
$\mathcal{A}$ we get that $[X(t)-X(s)]e\in\mathcal{A}$. In
particular, if $X(s_0)\in\mathcal{A}$ for some $s_0\in[a,b]$ then
right Segal's uniform continuity means that for each
$\varepsilon>0$ there is a projection $e\in\mathcal{A}\text{ with
}\tau(e)>1-\varepsilon$ such that the process $(X(t)e\colon
t\in[a,b])\subset\mathcal{A}$ is uniformly continuous in
$\|\cdot\|_{\infty}$-norm. The same holds of course for left
Segal's uniform continuity.

It seems worthwhile to say a few words about the terminology. The
term `Segal's convergence' was introduced by E.C. Lance in
\cite{La} in honour of I. Segal who first considered this mode of
convergence in his celebrated paper \cite{Se}. This notion
consists in the following: $x_n\to x$ \emph{in Segal's sense} if
for each $\varepsilon>0$ there is a projection
$e\in\mathcal{A}\text{ with }\tau(e^{\bot})<\varepsilon$ such that
$(x_n-x)e\in\mathcal{A}$ for sufficiently large $n,\text{ and
}\|(x_n-x)e\|_{\infty}\to 0$. In the definition above it is
assumed that $\tau$ is a faithful normal \emph{semifinite} trace
on $\mathcal{A}$. If $\tau$ is finite (as in our case) then
Segal's convergence becomes the so-called \emph{almost uniform
convergence} (which in the commutative case is \emph{via} Egorov's
theorem equivalent to convergence almost everywhere). Now the
similarity between Segal's (or in other words: almost uniform)
convergence and Segal's continuity is obvious and goes
(essentially) like that: in Segal's convergence we can find an
`arbitrarily large' projection $e$ such that $x_ne\to xe\text{
 in }\|\cdot\|_{\infty}$-norm, while in Segal's
continuity we can find an `arbitrarily large' projection $e$ such
that the process $(X(t)e\colon t\in[a,b])$ is uniformly continuous
in $\|\cdot\|_{\infty}$-norm. Accordingly, Segal's continuity
might also be called \emph{almost uniform continuity}.
\begin{remark}
It takes little effort to show that Segal's uniform continuity can
be given the following, equivalent but technically simpler, form:
\[\label{S-c}
\begin{aligned}
 &\begin{aligned}
 &\text{for each }\varepsilon>0\text{ there are a projection } e
 \in\mathcal{A}\\&\text{with }\tau(e)>1-\varepsilon,\text{and
 }\delta>0,
 \text{ such that for any }s,t\in[a,b]\\
 &\text{with }|t-s|<\delta,\text{ we have}
 \end{aligned}\notag\\
 &\begin{aligned}
 &[X(t)-X(s)]e\in\mathcal{A},& e[X(t)-X(s)]\in\mathcal{A}&\\
 &\text{and}\\
 &\|[X(t)-X(s)]e\|_{\infty}\leqslant\varepsilon,\qquad&
 \|e[X(t)-X(s)]\|_{\infty}\leqslant\varepsilon&.
 \end{aligned}
\end{aligned}\tag{S-cont}
\]
\end{remark}
Consider now a process $(X(t)\colon t\in[0,+\infty))$. It is
easily seen that the trajectories of this process are continuous
with probability one if and only if for each bounded interval
$[a,b]$ contained in $[0,+\infty)$ the trajectories of the process
$(X(t)\colon t\in[a,b])$ are uniformly continuous. In accordance
with the above observation we adopt the following definition.
\begin{definition}
 Let $(X(t)\colon t\in[0,+\infty))$ be a noncommutative stochastic
 process. We say that it is \emph{continuous in Segal's
 sense} if for any subinterval $[a,b]$ of the interval
 $[0,+\infty)$ the process $(X(t)\colon t\in[a,b])$ is uniformly continuous in
 Segal's sense, i.e. condition \eqref{S-c} is satisfied.
\end{definition}
The considerations above lead to one more notion of continuity.
Namely, the projection $e$ occurring in the definitions of left
and right Segal's continuity can be put on both sides.
Accordingly, we have
\begin{definition}
Let $(X(t)\colon t\in[0,+\infty))$ be a noncommutative stochastic
process. We say that it is \emph{weakly continuous in Segal's
sense} if for any subinteval $[a,b]$ of the interval $[0,+\infty)$
the process $(X(t)\colon t\in[a,b])$ is \emph{weakly uniformly
continuous in Segal's sense}, i.e. for each $\varepsilon>0$ there
is a projection $e$ in $\mathcal{A}\text{ with
}\tau(e)>1-\varepsilon$, having the property: for each $\eta>0$
there is $\delta>0$ such that for any $s,t\in[a,b]\text{ with
}|t-s|<\delta$, we have
\begin{equation*}
 e[X(t)-X(s)]e\in\mathcal{A}\qquad
 \text{and}\qquad\|e[X(t)-X(s)]e\|_{\infty}\leqslant\eta,
\end{equation*}
or equivalently, for each $\varepsilon>0$ there are a projection
$e\in\mathcal{A}$ with $\tau(e)>1-\varepsilon,\text{ and
 }\delta>0$, such that for any $s,t\in[a,b]
 \text{ with }|t-s|<\delta$, we have
\[
 e[X(t)-X(s)]e\in\mathcal{A}\qquad\text{and}\qquad
 \|e[X(t)-X(s)]e\|_{\infty}\leqslant\varepsilon.
\]
\end{definition}
It is clear that both left and right Segal's uniform continuity
imply weak Segal's uniform continuity; moreover, if $(X(t))$ is
left and $(Y(t))$ is right uniformly continuous in Segal's sense
then $(X(t)+Y(t))$ is weakly uniformly continuous in Segal's
sense. Obviously, in the commutative case all three modes of
continuity are equivalent.

\section{Stochastic integral}
In this section we shall prove the following
\begin{theorem}\label{Int}
Let $(X(t)\colon t\in[0,+\infty))$ be an $L^p$-valued
martingale,\linebreak $p>2$, and let $(f(t)\colon
t\in[0,+\infty))$ be an $\mathcal{A}$-valued adapted norm-bounded
on each interval $[a,b]$ continuous in Segal's sense process. Then
for each $t>0$ there exist stochastic integrals
\[
 Y(t)=\int_0^t dX(u)\,f(u)\qquad\text{and}\qquad Z(t)=\int_0^t f(u)\,dX(u)
\]
as elements of $L^2(\mathcal{A})$. Moreover, $(Y(t))\text{ and
}(Z(t))$ are martingales, and if $(X(t))$ is continuous either in
Segal's sense or in $\|\cdot\|_2$-norm then these martingales are
$L^2$-continuous.
\end{theorem}
\begin{proof}
Restrict attention to the left integral. Let \linebreak
$\theta(t)=\{0=t_0<t_1<\dots<t_m=t\}$ be a partition of $[0,t]$,
and let
\[
 S_{\theta(t)}^l(t)=\sum_{k=1}^m\Delta X(t_k)f(t_{k-1}).
\]
We shall show that for each $a>0$
\[
 \lim_{\|\theta(t)\|\to 0}S_{\theta(t)}^l(t)=Y(t)\qquad \text{
 uniformly in }t\in[0,a],
\]
that is for each $\varepsilon>0$ there is $\delta>0$ such that for
any $t\in[0,a]$ and any $\theta(t)\text{ with }
\|\theta(t)\|<\delta$ we have
\[
 \|S_{\theta(t)}^l(t)-Y(t)\|_2\leqslant\varepsilon.
\]
Put
\begin{align}
 M&=\sup_{0\leqslant u\leqslant a}\|f(u)\|_{\infty},\label{o1}\\
 K&=\tau(|X(a)|^2)-\tau(|X(0)|^2)\label{o2}\\
 q&=\frac{p}{p-2}\quad\text{so
 that}\quad\frac{1}{q}+\frac{2}{p}=1,\label{q}
\end{align}
and let $\alpha_p$ be as in \eqref{PX}. Take an arbitrary
$t\in[0,a]$, and let $\varepsilon>0$ be given. On account of
Segal's continuity of $f$ we can find a projection $e$ in
$\mathcal{A}$ with
\begin{align}\label{o3}
 &\tau(e)>1-\left[\frac{\varepsilon^2}{32M^2(\alpha_p\|X(a)\|_p)^2}\right]^q
 \qquad\text{i.e.}\notag\\
 &\tau(e^{\bot})<\left[\frac{\varepsilon^2}
 {32M^2(\alpha_p\|X(a)\|_p)^2}\right]^q,
\end{align}
and $\delta>0$ such that for any $t',t''\in[0,t]$ with
$|t'-t''|<\delta$, we have
\begin{equation}\label{e2}
 \begin{aligned}
  &\|e[f(t')-f(t'')]\|_{\infty}=\|[f(t')-f(t'')]^*e\|_{\infty}
  \leqslant\frac{\varepsilon^2}{32MK}\\
  &\|[f(t')-f(t'')]e\|_{\infty}
  \leqslant\frac{\varepsilon^2}{32MK}.
 \end{aligned}
\end{equation}
Let $\theta'(t)=\{0=t_0<t_1<\dots<t_m=t\}$ be an arbitrary
partition of $[0,t]$ with $\|\theta'(t)\|<\delta$ and let
$\theta''(t)$ be a partition of $[0,t]$ finer than $\theta'(t)$ .
Denote by $t_0^{(k)}, t_1^{(k)},\dots,t_{l_k}^{(k)}$ the points of
$\theta''(t)$ lying between $t_{k-1}\text{ and }t_k$, such that
$t_{k-1}=t_0^{(k)}<t_1^{(k)}<\dots<t_{l_k}^{(k)}=t_k$. We then
have
\begin{align*}
 S_{\theta''(t)}^l&=\sum_{k=1}^m\sum_{i=1}^{l_k}\Delta X(t_i^{(k)})
 f(t_{i-1}^{(k)})\\
 S_{\theta'(t)}^l&=\sum_{k=1}^m\Delta X(t_k)f(t_{k-1})=
 \sum_{k=1}^m\sum_{i=1}^{l_k}\Delta X(t_i^{(k)})f(t_{k-1}),
\end{align*}
so that
\[
 S_{\theta''(t)}^l-S_{\theta'(t)}^l=\sum_{k=1}^m\sum_{i=1}^{l_k}
 \Delta X(t_i^{(k)})[f(t_{i-1}^{(k)})-f(t_{k-1})].
\]
Consequently,
\begin{align*}
 \|S_{\theta''(t)}^l-S_{\theta'(t)}^l\|_2^2=&\tau\Big(\sum_{r=1}^m\sum_{j=1}^{l_r}
 [f(t_{j-1}^{(r)})-f(t_{r-1})]^*\Delta X(t_j^{(r)})^*\cdot\\
 &\cdot\sum_{k=1}^m\sum_{i=1}^{l_k}\Delta
 X(t_i^{(k)})[f(t_{i-1}^{(k)})-f(t_{k-1})]\Big)\\
 =&\tau\Big(\sum_{r=1}^m\sum_{j=1}^{l_r}\sum_{k=1}^m\sum_{i=1}^{l_k}
 [f(t_{j-1}^{(r)})-f(t_{r-1})]^*\Delta X(t_j^{(r)})^*\cdot\\
 &\cdot\Delta X(t_i^{(k)})[f(t_{i-1}^{(k)})-f(t_{k-1})]\Big)\\
 =&\sum_{r=1}^m\sum_{j=1}^{l_r}\sum_{k=1}^m\sum_{i=1}^{l_k}
 \tau([f(t_{i-1}^{(k)})-f(t_{k-1})]\cdot\\
 &\cdot[f(t_{j-1}^{(r)})-f(t_{r-1})]^*\Delta X(t_j^{(r)})^*\Delta X(t_i^{(k)})).
\end{align*}
For $k<r$ we have $t_{k-1}\leqslant
t_{i-1}^{(k)}<t_i^{(k)}\leqslant t_{r-1}\leqslant
t_{j-1}^{(r)}<t_j^{(r)}$, and the
\mbox{$\mathbb{E}_t$-in}\-variance of $\tau$ together with the
fact that $f(t_{i-1}^{(k)})-f(t_{k-1})$,
\mbox{$[f(t_{j-1}^{(r)})-f(t_{r-1})]^*$} and $\Delta X(t_i^{(k)})$
belong to $L^p(\mathcal{A}_{t_{j-1}^{(r)}})$ yield
\begin{align*}
 &\tau([f(t_{i-1}^{(k)})-f(t_{k-1})][f(t_{j-1}^{(r)})-f(t_{r-1})]^*
 \Delta X(t_j^{(r)})^*\Delta X(t_i^{(k)}))\\
 =&\tau(\mathbb{E}_{t_{j-1}^{(r)}}[f(t_{i-1}^{(k)})-f(t_{k-1})]
 [f(t_{j-1}^{(r)})-f(t_{r-1})]^*
 \Delta X(t_j^{(r)})^*\Delta X(t_i^{(k)}))\\
 =&\tau([f(t_{i-1}^{(k)})-f(t_{k-1})][f(t_{j-1}^{(r)})-f(t_{r-1})]^*
 [\mathbb{E}_{t_{j-1}^{(r)}}\Delta X(t_j^{(r)})]^*\Delta
 X(t_i^{(k)}))=0,
\end{align*}
since
\[
 \mathbb{E}_{t_{j-1}^{(r)}}\Delta X(t_j^{(r)})=\mathbb{E}_{t_{j-1}^{(r)}}
 (X(t_j^{(r)})-X(t_{j-1}^{(r)}))=\mathbb{E}_{t_{j-1}^{(r)}}X(t_j^{(r)})-
 X(t_{j-1}^{(r)})=0
\]
by martingale property. Analogously for $k>r$, thus we are left
only with the case $k=r$. For $i<j$ we have $t_{k-1}\leqslant
t_{i-1}^{(k)}<t_i^{(k)}\leqslant t_{j-1}^{(k)}<t_j^{(k)}$, and in
a similar fashion as above we obtain
\begin{align*}
 &\tau([f(t_{i-1}^{(k)})-f(t_{k-1})][f(t_{j-1}^{(k)})-f(t_{k-1})]^*
 \Delta X(t_j^{(k)})^*\Delta X(t_i^{(k)}))\\
 =&\tau(\mathbb{E}_{t_{j-1}^{(k)}}[f(t_{i-1}^{(k)})-f(t_{k-1})]
 [f(t_{j-1}^{(k)})-f(t_{k-1})]^*
 \Delta X(t_j^{(k)})^*\Delta X(t_i^{(k)}))\\
 =&\tau([f(t_{i-1}^{(k)})-f(t_{k-1})][f(t_{j-1}^{(k)})-f(t_{k-1})]^*
 [\mathbb{E}_{t_{j-1}^{(k)}}\Delta X(t_j^{(k)})]^*\Delta
 X(t_i^{(k)}))=0.
\end{align*}
The same goes for $i>j$, so finally we get
\begin{align*}
 \|S_{\theta''(t)}^l-S_{\theta'(t)}^l\|_2^2=&\tau\Big(\sum_{k=1}^m\sum_{i=1}^{l_k}
 [f(t_{i-1}^{(k)})-f(t_{k-1})][f(t_{i-1}^{(k)})-f(t_{k-1})]^*
 \cdot\\&\cdot\Delta X(t_i^{(k)})^*\Delta
 X(t_i^{(k)})\Big)=I_1+I_2,
\end{align*}
where
\begin{align*}
 I_1&=\tau\Big(\sum_{k=1}^m\sum_{i=1}^{l_k}
 e[f(t_{i-1}^{(k)})-f(t_{k-1})][f(t_{i-1}^{(k)})-f(t_{k-1})]^*
 |\Delta X(t_i^{(k)}|^2\Big),\\
 I_2&=\tau\Big(\sum_{k=1}^m\sum_{i=1}^{l_k}
 e^{\bot}[f(t_{i-1}^{(k)})-f(t_{k-1})][f(t_{i-1}^{(k)})-f(t_{k-1})]^*
 |\Delta X(t_i^{(k)}|^2\Big).
\end{align*}
For $I_1$ we have using \eqref{o1}, \eqref{o2} and \eqref{e2}
together with \eqref{e1}
\begin{align*}
 |I_1|&\leqslant\sum_{k=1}^m\sum_{i=1}^{l_k}
 |\tau(e[f(t_{i-1}^{(k)})-f(t_{k-1})][f(t_{i-1}^{(k)})-f(t_{k-1})]^*
 |\Delta X(t_i^{(k)})|^2)|\\
 &\leqslant\sum_{k=1}^m\sum_{i=1}^{l_k}\|e[f(t_{i-1}^{(k)})-f(t_{k-1})]
 \|_{\infty}\|[f(t_{i-1}^{(k)})-f(t_{k-1})]^*\|_{\infty}\cdot\\
 &\phantom{\sum\sum}\quad\cdot\||\Delta X(t_i^{(k)})|^2\|_1
 \leqslant\frac{\varepsilon^2}{32MK}2M\sum_{k=1}^m\sum_{i=1}^{l_k}
 \||\Delta X(t_i^{(k)})|^2\|_1 \\ &\leqslant\frac{\varepsilon^2}{16K}\tau
 (|X(t)|^2-|X(0)|^2)\leqslant\frac{\varepsilon^2}{16}.
\end{align*}
Divide $I_2$ into two parts: $I_2=I_2'+I_2''$, where
\begin{align*}
 I_2'&=\tau\Big(\sum_{k=1}^m\sum_{i=1}^{l_k}
 e^{\bot}[f(t_{i-1}^{(k)})-f(t_{k-1})][f(t_{i-1}^{(k)})-f(t_{k-1})]^*
 e^{\bot}|\Delta X(t_i^{(k)})|^2\Big),\\
 I_2''&=\tau\Big(\sum_{k=1}^m\sum_{i=1}^{l_k}
 e^{\bot}[f(t_{i-1}^{(k)})-f(t_{k-1})][f(t_{i-1}^{(k)})-f(t_{k-1})]^*
 e^|\Delta X(t_i^{(k)})|^2\Big).
\end{align*}
Then we have by \eqref{o1}
\begin{align*}
 |I_2'|&\leqslant\sum_{k=1}^m\sum_{i=1}^{l_k}
 |\tau(e^{\bot}[f(t_{i-1}^{(k)})-f(t_{k-1})][f(t_{i-1}^{(k)})-f(t_{k-1})]^*
 e^{\bot}|\Delta X(t_i^{(k)})|^2)|\\
 &=\sum_{k=1}^m\sum_{i=1}^{l_k}
 |\tau([f(t_{i-1}^{(k)})-f(t_{k-1})][f(t_{i-1}^{(k)})-f(t_{k-1})]^*
 e^{\bot}|\Delta X(t_i^{(k)})|^2e^{\bot})|\\
 &\leqslant\sum_{k=1}^m\sum_{i=1}^{l_k}
 \|([f(t_{i-1}^{(k)})-f(t_{k-1})][f(t_{i-1}^{(k)})-f(t_{k-1})]^*\|_{\infty}\cdot
 \\&\phantom{\leqslant\sum\sum}\cdot\|e^{\bot}|\Delta X(t_i^{(k)})|^2e^{\bot})\|_1\\
 &\leqslant4M^2\sum_{k=1}^m\sum_{i=1}^{l_k}\tau(e^{\bot}|\Delta
 X(t_i^{(k)})|^2e^{\bot})=4M^2\tau\Big(e^{\bot}\sum_{k=1}^m\sum_{i=1}^{l_k}|\Delta
 X(t_i^{(k)})|^2\Big)\\&\leqslant4M^2\Big\|e^{\bot}\sum_{k=1}^m\sum_{i=1}^{l_k}|\Delta
 X(t_i^{(k)})|^2)\Big\|_1.
\end{align*}
From H\"{o}lder's inequality applied to the last term above we
get, with $q$ given by \eqref{q},
\begin{align*}
 |I_2'|&\leqslant4M^2\Big\|e^{\bot}\Big\|_q\Big\|\sum_{k=1}^m\sum_{i=1}^{l_k}|\Delta
 X(t_i^{(k)})|^2\Big\|_{p/2}\\
 &=4M^2[\tau(e^{\bot})]^{1/q}\Big\|(\sum_{k=1}^m\sum_{i=1}^{l_k}|\Delta
 X(t_i^{(k)})|^2)^{1/2}\Big\|_p^2,
\end{align*}
and taking into account inequalities \eqref{PX}, \eqref{o3} we
finally obtain
\begin{align*}
 |I_2'|&\leqslant4M^2\frac{\varepsilon^2}{32M^2(\alpha_p\|X(a)\|_p)^2}(\alpha_p\|X(t)\|_p)^2
 \leqslant\frac{\varepsilon^2}{8}.
\end{align*}
For $I_2''$ we have using again \eqref{o1}, \eqref{e2} and
\eqref{o2} together with \eqref{e1}
\begin{align*}
 |I_2''|&\leqslant\sum_{k=1}^m\sum_{i=1}^{l_k}
 |\tau(e^{\bot}[f(t_{i-1}^{(k)})-f(t_{k-1})][f(t_{i-1}^{(k)})-f(t_{k-1})]^*
 e|\Delta X(t_i^{(k)})|^2)|\\
 &\leqslant\sum_{k=1}^m\sum_{i=1}^{l_k}
 \|e^{\bot}[f(t_{i-1}^{(k)})-f(t_{k-1})][f(t_{i-1}^{(k)})-f(t_{k-1})]^*e\|_{\infty}\cdot\\
 &\phantom{\leqslant\sum\sum}\cdot\||\Delta X(t_i^{(k)})|^2)\|_1\\
 &\leqslant\sum_{k=1}^m\sum_{i=1}^{l_k}2M\frac{\varepsilon^2}{32MK}\||\Delta
 X(t_i^{(k)})|^2\|_1=\frac{\varepsilon^2}{16K}\sum_{k=1}^m\sum_{i=1}^{l_k}
 \|\Delta
 X(t_i^{(k)})|^2\|_1\\&=\frac{\varepsilon^2}{16K}\tau(|X(t)|^2-|X(0)|^2)\leqslant
 \frac{\varepsilon^2}{16}.
\end{align*}
Thus
\begin{equation*}
 |I_2|\leqslant|I_2'|+|I_2''|\leqslant\frac{\varepsilon^2}{8}+\frac{\varepsilon^2}{16},
\end{equation*}
and consequently,
\begin{equation}\label{e3}
 \|S_{\theta''(t)}^l-S_{\theta'(t)}^l\|_2^2=I_1+I_2\leqslant\frac{\varepsilon^2}{16}+
 \frac{\varepsilon^2}{8}+\frac{\varepsilon^2}{16}=\frac{\varepsilon^2}{4}.
\end{equation}
Let now $\theta_1(t)\text{ and }\theta_2(t)$ be arbitrary
partitions
of $[0,t]$ such that\\
$\|\theta_1(t)\|<\delta,\,\|\theta_2(t)\|<\delta$, and let
$\theta''(t)=\theta_1(t)\cup\theta_2(t)$. Then we have by
\eqref{e3}
\[
 \|S_{\theta''(t)}^l(t)-S_{\theta_1(t)}^l(t)\|_2\leqslant
 \frac{\varepsilon}{2}\qquad\text{and}\qquad
 \|S_{\theta''(t)}^l(t)-S_{\theta_2(t)}^l(t)\|_2
 \leqslant\frac{\varepsilon}{2},
\]
so
\[
 \|S_{\theta_1(t)}^l(t)-S_{\theta_2(t)}^l(t)\|_2\leqslant\varepsilon,
\]
which means that the net $\{S_{\theta(t)}^l(t)\}$ satisfies the
Cauchy condition\\as $\|\theta(t)\|\to 0$ uniformly in
$t\in[0,a]$, proving the existence of the left integral
$Y(t)=\int_0^t dX(t)\,f(t)$, together with the uniform convergence
of $S_{\theta(t)}(t)\text{ to }Y(t) \text{ in } \|\cdot\|_2$-norm.
The existence of the right integral $Z(t)=\int_0^t f(t)\,dX(t)$ is
proved in virtually the same way.

Now we shall show that $(Y(t)\colon t\in[0,+\infty))$ is a
martingale. Fix $t>0$ and take an arbitrary $s<t$ ($t$ being
fixed, so we suppress in our notation the dependence of
$\theta\text{ and }S_{\theta}^l \text{ on }t$). We have
\[
 \int_0^t dX(u)\,f(u)=\lim_{\|\theta\|\to 0}S_{\theta}^l.
\]
We may assume that $s$ is one of the points of each partition
\linebreak $\theta=\{0=t_0<t_1<\cdots<t_m=t\}$, say $s=t_k$. Then
we have
\begin{align*}
 &\mathbb{E}_s S_{\theta}^l=\mathbb{E}_s(\sum_{i=1}^k
 [X(t_i)-X(t_{i-1})]f(t_{i-1})+\sum_{i=k+1}^m
 [X(t_i)-X(t_{i-1})]f(t_{i-1}))\\
 &=\sum_{i=1}^k\mathbb{E}_s[X(t_i)-X(t_{i-1})]f(t_{i-1})+
 \sum_{i=k+1}^m\mathbb{E}_s[X(t_i)-X(t_{i-1})]f(t_{i-1}).
\end{align*}
For $i\leqslant k$ we have $t_i\leqslant s$, and thus
\[
 \mathbb{E}_s[X(t_i)-X(t_{i-1})]f(t_{i-1})=[X(t_i)-X(t_{i-1})]f(t_{i-1}),
\]
while for $i>k$ we have $t_{i-1}\geqslant s$, and thus
\begin{align*}
 &\mathbb{E}_s[X(t_i)-X(t_{i-1})]f(t_{i-1})=
 \mathbb{E}_s\mathbb{E}_{t_{i-1}}[X(t_i)-X(t_{i-1})]f(t_{i-1})\\
 =&\mathbb{E}_s(\mathbb{E}_{t_{i-1}}[X(t_i)-X(t_{i-1})])f(t_{i-1})=0
\end{align*}
by martingale property. Consequently,
\begin{equation}\label{e4}
 \mathbb{E}_s
 S_{\theta}^l=\sum_{i=1}^k[X(t_i)-X(t_{i-1})]f(t_{i-1}).
\end{equation}
But the sum on the right hand side of \eqref{e4} is an integral
sum for the integral $\int_0^s dX(u)\,f(u)$, and passing to the
limit in \eqref{e4} yields
\[
 \mathbb{E}_s\int_0^t dX(u)\,f(u)=\int_0^s dX(u)\,f(u),
\]
which shows that $(Y(t))$ is a martingale. Analogously for
$(Z(t))$.

Now we shall prove the $\|\cdot\|_2$-continuity of $(Y(t))$ in an
arbitrary interval $[0,a]$. Let
$\theta_n=\{0=t_1^{(n)}<t_2^{(n)}<\dots<t_{m_n}^{(n)}=a\}$ be a
sequence of partitions of $[0,a]$ such that
$\theta_n\subset\theta_{n+1}\text{ and }\|\theta_n\|\to 0$. For an
arbitrary $t\in[0,a]$ put
$\theta_n(t)=(\theta_n\cap[0,t])\cup\{t\}$, and
\[
 S_n(t)=\sum_{0\leqslant t_{k-1}^{(n)}<t_k^{(n)}\leqslant
 t}\Delta X(t_k^{(n)})f(t_{k-1}^{(n)})+[X(t)-X(t'_n)]f(t'_n),
\]
where $t'_n=\max\{t_k^{(n)}\colon t_k^{(n)}\leqslant t\}$. Then
\[
 S_n(t)=S_{\theta_n(t)}^l(t)
\]
in the notation from the first part of the proof, and
\[
 \lim_{n\to\infty}S_n(t)=Y(t)\qquad \text{ uniformly in }
 t\in[0,a].
\]

First we show the $\|\cdot\|_2$-continuity of $S_n$. Take
arbitrary $s,t\in[0,a]$, $s<t$, such that $|t-s|<\|\theta_n\|$. We
have three possibilities: \par (i)
$s\in\theta_n,\;t\notin\theta_n$.
\\Then $t'_n=s$, and
\[
 S_n(t)-S_n(s)=[X(t)-X(s)]f(s).
\]
\par
(ii) $s\notin\theta_n,\;t\in\theta_n$.
\\Then $t'_n=t$, and
\begin{align*}
 S_n(t)-S_n(s)&=[X(t)-X(s'_n)]f(s'_n)-[X(s)-X(s'_n)]f(s'_n)\\
 &=[X(t)-X(s)]f(s'_n),
\end{align*}
where $s'_n=\max\{t_k^{(n)}\colon t_k^{(n)}\leqslant s\}$.
\par
(iii) $s\notin\theta_n,\;t\notin\theta_n$.
\\Then $s'_n<s<t'_n<t$, and
\begin{equation}\label{Sc}
\begin{aligned}
 &S_n(t)-S_n(s)=[X(t)-X(t'_n)]f(t'_n)\\+&[X(t'_n)-X(s'_n)]f(s'_n)
 -[X(s)-X(s'_n)]f(s'_n)\\=&[X(t)-X(t'_n)]f(t'_n)+[X(t'_n)-X(s)]f(s'_n).
\end{aligned}
\end{equation}

Note that case (iii) contains the other two, so we have in general
\begin{equation}\label{uc}
 \|S_n(t)-S_n(s)\|_2\leqslant
 M\|X(t)-X(t'_n)\|_2+M\|X(t'_n)-X(s)\|_2.
\end{equation}
Consequently, we only need to estimate the expression
$\|X(u)-X(v)\|_2$ for $s\leqslant v<u\leqslant t$.

If $(X(t))$ is $\|\cdot\|_2$-continuous then obviously this can be
made arbitrarily small for $s,t$ sufficiently close to each other.

If $(X(t))$ is continuous in Segal's sense then for each
$\varepsilon>0$ we can find a projection $e$ in $\mathcal{A}$ with
\[
 \tau(e^{\bot})<\left[\frac{\varepsilon}{8M\|X(a)\|_p}\right]^{q'},
 \qquad\text{where}\qquad \frac{1}{q'}+\frac{1}{p}=\frac{1}{2},
\]
and $\delta>0$ such that for $|t'-t''|<\delta$ we have
\[
 \|e[X(t')-X(t'')]\|_{\infty}<\frac{\varepsilon}{4M}.
\]
If $|t-s|<\delta$ then also $|u-v|<\delta$, and
\begin{align*}
 \|X(u)-X(v)\|_2&\leqslant\|e[X(u)-X(v)]\|_2+
 \|e^{\bot}[X(u)-X(v)]\|_2 \\ &\leqslant
 \|e[X(u)-X(v)]\|_{\infty}+ \|e^{\bot}[X(u)-X(v)]\|_2\\&\leqslant
 \frac{\varepsilon}{4M}+ \|e^{\bot}[X(u)-X(v)]\|_2.
\end{align*}
From H\"{o}lder's inequality we obtain
\begin{align*}
 &\|e^{\bot}[X(u)-X(v)]\|_2\leqslant\|e^{\bot}\|_{q'}\|X(u)-X(v)\|_p\\
 &\leqslant
 2\|X(a)\|_p\left[\tau(e^{\bot})\right]^{1/q'}<\frac{\varepsilon}{4M},
\end{align*}
thus finally
\[
 \|X(u)-X(v)\|_2<\frac{\varepsilon}{2M}.
\]
This estimation shows that for $|t-s|<\delta\wedge\|\theta_n\|$
the inequality in \eqref{uc} takes the form
\[
 \|S_n(t)-S_n(s)\|_2<M\frac{\varepsilon}{2M}+M\frac{\varepsilon}{2M}=\varepsilon,
\]
showing the uniform $\|\cdot\|_2$-continuity of $S_n$. Since
$S_n\rightrightarrows Y$ in $[0,a]$ we obtain the
$\|\cdot\|_2$-continuity of $Y$.
\end{proof}
\begin{remark}
It is seen from the above proof that for the existence of the left
integral it suffices that $f$ be `left Segal's continuous' while
for the existence of the right integral it suffices that $f$ be
`right Segal's continuous'.
\end{remark}
Our next aim is to show a noncommutative counterpart of the known
classical result saying that a stochastic integral with the
integrator being a continuous martingale is continuous. To this
end, we begin with a result which may be looked upon as a
noncommutative generalisation of one of the classical martingale
inequalities. Our attention will be restricted to its simplest
version for a finite martingale, which suffices for the purposes
of this paper; however, it is worth mentioning that a result of
this type can be obtained also in a more general setting. The idea
of the proof is an adaptation of the classical method to the
noncommutative setup, and has already been used by C.J.K. Batty in
a slightly different context for proving `noncommutative
Kolmogorov's inequality' for sums of `independent noncommutative
random variables' (cf. \cite[Proposition 5.1]{Ba}).
\begin{proposition}\label{Ki}
Let $(X_1,\dots,X_m)$ be an $L^2$-martingale. Then for each
$\varepsilon>0$ there is a projection $e\in\mathcal{A}$ such that
\[
 \tau(e^{\bot})<\frac{\|X_m\|_2^2}{\varepsilon^2},
\]
and
\[
 \|eX_n\|_{\infty}\leqslant\varepsilon \qquad\text{\em for each
}\quad n=1,\dots,m.
\]
The same conclusion holds also in the `right version' with the
projection $e$ put to the right of $X_n$.
\end{proposition}
\begin{proof}
We shall prove the `left version', the proof of the other one
being \emph{mutatis mutandis} the same. Let us start with some
simple remarks. Let $x\in L^1(\mathcal{A})$ be a positive operator
with its spectral decomposition
\[
 x=\int_0^{\infty}\lambda\,e(d\lambda).
\]
For each $\eta>0$ we have
\begin{equation}\label{sp}
\begin{aligned}
 x&=\int_{[0,\eta)}\lambda\,e(d\lambda)+\int_{[\eta,+\infty)}\lambda\,
 e(d\lambda)\\ &\geqslant\int_{[\eta,+\infty)}\lambda\,e(d\lambda)
 \geqslant\eta e([\eta,+\infty)),
\end{aligned}
\end{equation}
and
\[
 e([0,\eta))xe([0,\eta))=\int_{[0,\eta)}\lambda\,e(d\lambda).
\]
In particular, we obtain `noncommutative Chebyshev's inequality'
\[
 \tau(e([\eta,+\infty)))\leqslant\frac{\tau(x)}{\eta},
\]
and the estimation
\[
 \|e([0,\eta))xe([0,\eta))\|_{\infty}\leqslant\eta.
\]
Now let $f$ be an arbitrary projection in $\mathcal{A}$. Write the
spectral decomposition of $fxf$
\[
 fxf=\int_0^{\infty}\lambda\,h(d\lambda).
\]
Since $f$ commutes with $fxf$ it follows that $f$ commutes with
the spectral measure $h$ of $fxf$, consequently
\begin{equation}\label{spi}
 \|fh([0,\eta))xh([0,\eta))f\|_{\infty}=
 \|h([0,\eta))fxfh([0,\eta))\|_{\infty}\leqslant\eta.
\end{equation}
To facilitate notation, we agree to denote the value of the
spectral measure of the operator $x\geqslant 0$ on a Borel set
$Z\subset\mathbb{R}$ by $e_Z(x)$. Let an arbitrary $\varepsilon>0$
be given. Define inductively projections $e_n,\,f_n$ for
$n=1,\dots,m$, by
\begin{flalign*}
 e_1&=e_{[0,\varepsilon^2)}(|X_1^*|^2), &f_1&=e_1\\
 e_2&=e_{[0,\varepsilon^2)}(f_1|X_2^*|^2f_1),
 &f_2&=e_1\wedge e_2\\ \intertext{\dotfill}
 e_n&=e_{[0,\varepsilon^2)}(f_{n-1}|X_n^*|^2f_{n-1}),
 &f_n&=e_1\wedge\ldots\wedge e_n.
\end{flalign*}
We have $e_n\in\mathcal{A}_n$, and for $i<j$
\begin{equation}\label{pi}
 e_j^{\bot}\leqslant f_{j-1}\leqslant f_i\leqslant e_i,
\end{equation}
thus in particular $e_i^{\bot}e_j^{\bot}=0\text{ for }i\ne j$. Put
$e=f_m=e_1\wedge\ldots\wedge e_m$. Then
\[
 e^{\bot}=\bigvee_{n=1}^me_n^{\bot}=\sum_{n=1}^me_n^{\bot}.
\]
Inequality \eqref{sp} yields
\[
 \varepsilon^2e_n^{\bot}=\varepsilon^2e_{[\varepsilon^2,+\infty)}
 (f_{n-1}|X_n^*|^2f_{n-1})\leqslant f_{n-1}|X_n^*|^2f_{n-1},
\]
and thus multiplying both sides of the above inequality by
$e_n^{\bot}$ we obtain from \eqref{pi} and the fact that
$(|X_n^*|^2\colon n=1,\dots,m)$ is a submartingale
\[
 \varepsilon^2e_n^{\bot}\leqslant e_n^{\bot}|X_n^*|^2e_n^{\bot}
 \leqslant e_n^{\bot}\mathbb{E}_n|X_m^*|^2e_n^{\bot}=\mathbb{E}_n
 (e_n^{\bot}|X_m^*|^2e_n^{\bot}).
\]
Consequently,
\begin{align*}
 \varepsilon^2\tau(e_n^{\bot})\leqslant\tau(\mathbb{E}_n(e_n^{\bot}
 |X_m^*|^2e_n^{\bot}))=\tau(e_n^{\bot}|X_m^*|^2e_n^{\bot})
 =\tau(e_n^{\bot}|X_m^*|^2).
\end{align*}
Thus
\begin{align*}
 \varepsilon^2\tau(e^{\bot})&=\varepsilon^2\tau\Big(\bigvee_{n=1}^me_n^{\bot}\Big)
 =\varepsilon^2\sum_{n=1}^m\tau(e_n^{\bot})
 \leqslant\tau\Big(\Big(\sum_{n=1}^me_n^{\bot}\Big)|X_m^*|^2\Big)
 \\&=\tau(e^{\bot}|X_m^*|^2)\leqslant\tau(|X_m^*|^2)=\|X_m^*\|_2^2=\|X_m\|_2^2,
\end{align*}
which gives the desired estimation of $\tau(e^{\bot})$. For the
norm we have on account of \eqref{spi}
\begin{align*}
 &\|eX_n\|_{\infty}^2=\|X_n^*e\|_{\infty}^2=\|e|X_n^*|^2e\|_{\infty}^2\\
 =&\|f_{n-1}e|X_n^*|^2ef_{n-1}\|_{\infty}^2\leqslant
 \|f_{n-1}e_n|X_n^*|^2e_nf_{n-1}\|_{\infty}^2\\
 =&\|f_{n-1}e_{[0,\varepsilon^2)}(f_{n-1}|X_n^*|^2f_{n-1})
 |X_n^*|^2e_{[0,\varepsilon^2)}(f_{n-1}|X_n^*|^2f_{n-1})f_{n-1}\|_{\infty}^2
 \leqslant\varepsilon^2,
\end{align*}
which proves the claim.
\end{proof}
\begin{theorem}\label{Scm}
Let $(X(t))\text{ and }(f(t))$ be as in Theorem ~\ref{Int}. If
$(X(t))$ is left continuous in Segal's sense then $(Y(t))$ is left
continuous in Segal's sense. If $(X(t))$ is right continuous in
Segal's sense then $(Z(t))$ is right continuous in Segal's sense.
\end{theorem}
\begin{proof}
Again we restrict attention to the `left' case. We shall prove the
left uniform Segal's continuity of $(Y(t))$ in an arbitrary
interval $[0,a]$.\par Fix a positive integer $n$. It is easily
seen that $(S_n(t)\colon t\in[0,a])$ is a martingale; moreover
equality \eqref{Sc} shows that if $(X(t)\colon t\in[0,a])$ is
uniformly left continuous in Segal's sense then $(S_n(t)\colon
t\in[0,a])$ is also uniformly left continuous in Segal's
sense.\par Let now $m,\,n$ be arbitrary fixed positive integers.
The process $(S_n(t)-S_m(t)\colon t\in[0,a])$ is a uniformly left
continuous in Segal's sense martingale. For any given
$\varepsilon_{n,m}>0\text{ let } f_{n,m}$ be a projection in
$\mathcal{A}$ such that $\tau(f_{n,m}^{\bot})<\varepsilon_{n,m}$,
and the processes \\ $(f_{n,m}[S_n(t)-S_m(t)]\colon
t\in[0,a]),\,(f_{n,m}S_n(t)\colon t\in[0,a])$,\\
$(f_{n,m}S_m(t)\colon t\in[0,a])\subset\mathcal{A}$ are uniformly
continuous in $\|\cdot\|_{\infty}$-norm. In particular, there is
$\delta_{n,m}>0$ such that for all $t',t''\in[0,a]$ with
$|t'-t''|<\delta_{n,m}$ we have
\[
 \|f_{n,m}([S_n(t')-S_m(t'')]-[S_n(t'')-S_m(t'')])\|_{\infty}
 \leqslant\frac{\varepsilon_{n,m}}{2}.
\]
Choose points $0=t_0<t_1<\dots<t_r=a$ such that \linebreak
$\max_{1\leqslant i\leqslant r}(t_i-t_{i-1})<\delta$. For the
martingale \\ $(S_n(t_i)-S_m(t_i)\colon i=0,1,\dots,r)$ we infer
on account of Proposition ~\ref{Ki} that there exists a projection
$q_{n,m}\in\mathcal{A}$ with
\[
 \tau(q_{n,m}^{\bot})<\frac{4\|S_n(a)-S_m(a)\|_2^2}{\varepsilon_{n,m}^2}
\]
such that
\[
 \|q_{n,m}[S_n(t_i)-S_m(t_i)]\|_{\infty}\leqslant
 \frac{\varepsilon_{n,m}}{2},\qquad\text{for}\qquad i=0,1,\dots,r.
\]
Put $e_{n,m}=f_{n,m}\wedge q_{n,m}$. Then
\[
 \tau(e_{n,m})<\varepsilon_{n,m}+\frac{4\|S_n(a)-S_m(a)\|_2^2}
 {\varepsilon_{n,m}^2},
\]
and for each $t\in[0,a]$ there is $t_i$ such that
$|t-t_i|<\delta$, so
\begin{align*}
 &\|e_{n,m}[S_n(t)-S_m(t)]\|_{\infty}\leqslant\|e_{n,m}([S_n(t)-S_m(t)]
 -[S_n(t_i)-S_m(t_i)])\|_{\infty}\\
 +&\|e_{n,m}[S_n(t_i)-S_m(t_i)]\|_{\infty}\leqslant
 \|f_{n,m}([S_n(t)-S_m(t)]+\\-&[S_n(t_i)-S_m(t_i)])\|_{\infty}
 +\|q_{n,m}[S_n(t_i)-S_m(t_i)]\|_{\infty}\leqslant\frac{\varepsilon_{n,m}}{2}
 +\frac{\varepsilon_{n,m}}{2}=\varepsilon_{n,m}.
\end{align*}
Moreover, the processes $(e_{n,m}S_n(t)\colon t\in[0,a])$ and \\
$(e_{n,m}S_m(t)\colon t\in[0,a])$ are uniformly continuous in
$\|\cdot\|_{\infty}$-norm.

Let an arbitrary $\varepsilon>0$ be given. We have $S_n(a)\to
Y(a)$ in \linebreak $\|\cdot\|_2$-norm, so we can find a
subsequence $\{n_k\}$ such that
\[
 \|S_{n_{k+1}}(a)-S_{n_k}(a)\|_2<\frac{\varepsilon^3}{2^{3k+5}},\qquad
 k=1,2,\dots.
\]
Apply our previous considerations to the martingale\\
$(S_{n_{k+1}}(t)-S_{n_k}(t)\colon t\in[0,a])$. For
$\varepsilon_k=\varepsilon/2^{k+1}$ there is a projection
$e_k\in\mathcal{A}$ with
\[
 \tau(e_k^{\bot})<\frac{\varepsilon}{2^{k+1}}+\frac{4(2^{k+1})^2
 \|S_{n_{k+1}}(a)-S_{n_k}(a)\|_2^2}{\varepsilon^2}
 <\frac{\varepsilon}{2^{k+1}}+\frac{\varepsilon}{2^{k+1}}
 =\frac{\varepsilon}{2^k},
\]
such that for each $t\in[0,a]$
\[
 \|e_k[S_{n_{k+1}}(t)-S_{n_k}(t)]\|_{\infty}\leqslant
 \frac{\varepsilon}{2^{k+1}}.
\]
Put
\[
 e=\bigwedge_{k=1}^{\infty}e_k.
\]
Then $\tau(e^{\bot})<\varepsilon$, and for each $t\in[0,a]$
\begin{equation}\label{ucs}
 \|e[S_{n_{k+1}}(t)-S_{n_k}(t)]\|_{\infty}\leqslant\frac
 {\varepsilon}{2^{k+1}}.
\end{equation}
Since the processes $(e_kS_{n_k}(t)\colon
t\in[0,a]),\:k=1,2,\dots$ are uniformly continuous in
$\|\cdot\|_{\infty}$-norm it follows that the processes\\
$(eS_{n_k}(t)\colon t\in[0,a]),\:k=1,2,\dots$ are also uniformly
continuous in \linebreak $\|\cdot\|_{\infty}$-norm. Condition
\eqref{ucs} says that the sequence of processes \linebreak
$(eS_{n_k}(t)\colon t\in[0,a]),\:k=1,2,\dots$ is Cauchy in
$\|\cdot\|_{\infty}$-norm uniformly for $t\in[0,a]$. Since
\[
 eS_{n_k}(t)\to eY(t)\qquad\text{\em in }\|\cdot\|_2-\text{\em
 norm}
\]
it follows that
\[
 eS_{n_k}(t)\to eY(t)\qquad\text{\em in
 }\|\cdot\|_{\infty}-\text{\em norm}\quad\text{uniformly for
 }t\in[0,a],
\]
and the $\|\cdot\|_{\infty}$-norm continuity of
$(eS_{n_k}(t)\colon t\in[0,a])$ yields the norm continuity of
$(eY(t)\colon t\in[0,a])$ which proves the claim.
\end{proof}

\section{Quantum Doob-Meyer decomposition}
Let $(X(t)\colon t\in[0,+\infty))$ be an $L^2$-martingale. Then
the process $(|X(t)|^2\colon t\in[0,+\infty))$ is a submartingale,
and a Doob-Meyer decomposition is given by the representation
\begin{equation}\label{DM}
 |X(t)|^2=M(t)+A(t), \qquad t\in[0,+\infty),
\end{equation}
where $(M(t))$ is a martingale, and $(A(t))$ is an increasing
positive process, i.e. $0\leqslant A(s)\leqslant A(t)\text{ for
}0\leqslant s\leqslant t$. This decomposition has been obtained in
many concrete situations (cf.
\cite{BGW,BSW1,BSW2,BW1,BW2,BW3,BW4,BW5}); in general, a
sufficient condition for the existence of a Doob-Meyer
decomposition was given in \cite{BSW1} in the following form.\par
An $L^1$-process $(Y(t)\colon t\in[0,+\infty))$ is said to be of
\emph{class} $D$ if the set
\[
 \Big\{\sum_{k=1}^m \mathbb{E}_{t_{k-1}}(Y(t_k)-Y(t_{k-1}))\colon
 0\leqslant t_0<t_1<\dots<
 t_m,\;m=1,2,\dots\Big\}
\]
is weakly relatively compact. If $(|X(t)|^2)$ is of class $D$ then
it has a Doob-Meyer decomposition.\par We shall refer to this
condition in a slightly modified form; however, since in the proof
some specific features of the construction of the decomposition
will be exploited, it seems preferable to present it in full
detail. First we shall show the crucial property of the
submartingale $(|X(t)|^2)$, namely that it is of class $D$ (in our
case even with a stronger compactness requirement) on each finite
interval.
\begin{proposition}\label{wc}
Let $(X(t)\colon t\in[0,+\infty))$ be an $L^p$-martingale, $p>2$.
Then for each $a>0$ the set
\begin{align*}
 \mathcal{S}=\Big\{&\sum_{k=1}^m
 \mathbb{E}_{t_{k-1}}(|X(t_k)|^2-|X(t_{k-1})|^2)\colon 0\leqslant
 t_0<t_1<\dots<t_m=a,\\&m=1,2,\dots\Big\}\subset L^{p/2}(\mathcal{A})
\end{align*}
is weakly relatively compact.
\end{proposition}
\begin{proof}
Take arbitrary $0\leqslant t_0<t_1<\dots<t_m=a$. On account of
\eqref{e0} we have
\[
 \sum_{k=1}^m
 \mathbb{E}_{t_{k-1}}(|X(t_k)|^2-|X(t_{k-1})|^2)=\sum_{k=1}^m
 \mathbb{E}_{t_{k-1}}|\Delta X(t_k)|^2.
\]
The main result of \cite{J}, Theorem 0.1, says that there exists a
constant $c_{p/2}$ depending only on $p$ such that
\[
 \Big\|\sum_{k=1}^m
 \mathbb{E}_{t_{k-1}}\Delta|X(t_k)|^2\Big\|_{p/2}
 \leqslant c_{p/2}\Big\|\sum_{k=1}^m
 |\Delta X(t_k)|^2\Big\|_{p/2},
\]
consequently we obtain by \eqref{PX}
\begin{align*}
 &\Big\|\sum_{k=1}^m
 \mathbb{E}_{t_{k-1}}(|X(t_k)|^2-|X(t_{k-1})|^2)\Big\|_{p/2}=\Big\|\sum_{k=1}^m
 \mathbb{E}_{t_{k-1}}|\Delta X(t_k)|^2\Big\|_{p/2}\\ \leqslant& c_{p/2}\Big\|\sum_{k=1}^m
 |\Delta X(t_k)|^2\Big\|_{p/2}=c_{p/2}\Big[\tau\Big(\Big(\sum_{k=1}^m
 |\Delta X(t_k)|^2\Big)^{p/2}\Big)\Big]^{2/p}\\
 =&c_{p/2}\Big\|\Big(\sum_{k=1}^m|\Delta X(t_k)|^2\Big)^{1/2}\Big\|_p^2\leqslant
 c_{p/2}\,\alpha_p^2\|X(a)\|_p^2,
\end{align*}
which means that $\mathcal{S}$ is norm-bounded. Since
$L^{p/2}(\mathcal{A})$ is the dual space to $L^q(\mathcal{A})$
with $q$ given by \eqref{q}, and \emph{vice versa}, the conclusion
follows.
\end{proof}
Now we are ready to prove the existence of a Doob-Meyer
decomposition. The main idea of the proof is the same as in the
classical Rao's proof (cf. \cite[Theorem 4.10]{KS} or \cite{R} ),
however some additional refinements will be needed.
\begin{theorem}[Doob-Meyer decomposition]
Let $(X(t)\colon t\in[0,+\infty))$ be an $L^p$-martingale, $p>2$.
Then there exists a Doob-Meyer decomposition for $(|X(t)|^2\colon
t\in[0,+\infty))$.
\end{theorem}
\begin{proof}
Let $\mathbb{D}=\{t_k^{(n)}=\frac{k}{2^n}\colon k,n=1,2,\dots\}$
be the set of dyadic numbers. In the remaining part of the proof,
whenever the symbol $t_k^{(n)}$ is used, it will always be assumed
that $t_k^{(n)}\in\mathbb{D}$. Fix an arbitrary positive integer
$m$. For each $n=1,2,\dots$ and each dyadic number $u$ in
$[m-1,m]$ put
\[
 S_n^{(m)}(u)=\sum_{m-1\leqslant t_{k-1}^{(n)}<u}
 \mathbb{E}_{t_{k-1}^{(n)}}(|X(t_k^{(n)})|^2-|X(t_{k-1}^{(n)})|^2),
\]
in particular
\begin{equation*}
 S_n^{(m)}(m)=\sum_{k=(m-1)2^n+1}^{m2^n}
 \mathbb{E}_{t_{k-1}^{(n)}}(|X(t_k^{(n)})|^2-|X(t_{k-1}^{(n)})|^2).
\end{equation*}
By Proposition~\ref{wc} the sequence $\{S_n^{(m)}(m)\colon
n=1,2,\dots\}$ is weakly relatively compact. Let $A_m\in
L^{p/2}(\mathcal{A}_m)$ be a limit point of this sequence. By the
Eberlein-\u{S}mulian theorem there is a subsequence $\{n_r\}$
(depending on $m$) such that
\[
 \lim_{r\to\infty}S_{n_r}^{(m)}(m)=A_m \qquad \text{\em weakly}.
\]
From the weak continuity of conditional expectation we obtain
\begin{equation}\label{e5}
 \begin{aligned}
 &\;\mathbb{E}_{m-1}A_m\\
 =&\lim_{r\to\infty}\mathbb{E}_{m-1}\sum_{k=(m-1)2^{n_r}+1}^{m2^{n_r}}
 \mathbb{E}_{t_{k-1}^{(n_r)}}(|X(t_k^{(n_r)})|^2-|X(t_{k-1}^{(n_r)})|^2)\\
 =&\lim_{r\to\infty}\mathbb{E}_{m-1}(|X(m)|^2-|X(m-1)|^2)\\
 =&\;\mathbb{E}_{m-1}|X(m)|^2-|X(m-1)|^2.
 \end{aligned}
\end{equation}
Define a process $(A(t)\colon t\in[0,+\infty))$ as follows. For
$m=1,2,\dots$ and $t\in[0,+\infty)$ put
\begin{equation}\label{e6}
 A(t)=|X(t)|^2-\mathbb{E}_t|X(m)|^2+\mathbb{E}_t(A_1+\cdots+A_m),\qquad
 t\in[m-1,m].
\end{equation}
Verify first the correctness of this definition. Computing $A(m)$
in the interval $[m-1,m]$, we obtain
\[
 A(m)=|X(m)|^2-\mathbb{E}_m|X(m)|^2+\mathbb{E}_m(A_1+\cdots+A_m)=A_1+\cdots+A_m.
\]
For the interval $[m,m+1]$ we have
\begin{align*}
 A(m)&=|X(m)|^2-\mathbb{E}_m|X(m+1)|^2+\mathbb{E}_m(A_1+\cdots+A_{m+1})\\
 &=|X(m)|^2-\mathbb{E}_m|X(m+1)|^2+\mathbb{E}_mA_{m+1}+(A_1+\cdots+A_m)\\
 &=A_1+\cdots+A_m,
\end{align*}
because
\[
|X(m)|^2-\mathbb{E}_m|X(m+1)|^2+\mathbb{E}_mA_{m+1}=0
\]
by \eqref{e5}. Thus $(A(t))$ is well defined. Moreover, $(A(t))$
is right-continuous in $\|\cdot\|_{p/2}$-norm as a sum of the
submartingale $(|X(t)|^2)$ and two martingales, all being
$\|\cdot\|_{p/2}$-right-continu\-o\-us as noted in Section~1.
Putting $m=1,\,t=0$ in \eqref{e6} gives
\[
 A(0)=|X(0)|^2-\mathbb{E}_0|X(1)|^2+\mathbb{E}_0A_1,
\]
and putting $m=1$ in \eqref{e5} gives
\[
 \mathbb{E}_0A_1=\mathbb{E}_0|X(1)|^2-|X(0)|^2,
\]
which shows that $A(0)=0$.\par Let $u$ be an arbitrary dyadic
number in $[m-1,m]$. For sufficiently large $n$ (so large that
$u\in\{t_k^{(n)}\colon k=(m-1)2^n+1,\dots,m2^n\}$) the following
equality holds
\begin{align*}
 \mathbb{E}_uS_n^{(m)}(m)=&\sum_{m-1\leqslant t_{k-1}^{(n)}<u}
 \mathbb{E}_{t_{k-1}^{(n)}}(|X(t_k^{(n)})|^2-|X(t_{k-1}^{(n)})|^2)\\
 +&\sum_{u\leqslant t_{k-1}^{(n)}<m}
 \mathbb{E}_u(|X(t_k^{(n)})|^2-|X(t_{k-1}^{(n)})|^2)\\
 =&\sum_{m-1\leqslant t_{k-1}^{(n)}<u}\mathbb{E}_{t_{k-1}^{(n)}}
 (|X(t_k^{(n)})|^2-|X(t_{k-1}^{(n)}|^2))\\+&\;\mathbb{E}_u(|X(m)|^2-|X(u)|^2)\\
 =&S_n^{(m)}(u)+\mathbb{E}_u|X(m)|^2-|X(u)|^2.
\end{align*}
Passing to the limit along the sequence $\{n_r\}$ in the above
equality, we get
\[
 \mathbb{E}_uA_m=\lim_{r\to\infty}S_{n_r}^{(m)}(u)+\mathbb{E}_u|X(m)|^2-|X(u)|^2,
\]
and thus
\begin{equation}\label{e7}
 \begin{aligned}
 \lim_{r\to\infty}S_{n_r}^{(m)}(u)&=|X(u)|^2-\mathbb{E}_u|X(m)|^2+\mathbb{E}_uA_m\\
 &=A(u)-\mathbb{E}_u(A_1+\cdots+A_{m-1})\\
 &=A(u)-(A_1+\cdots+A_{m-1}).
 \end{aligned}
\end{equation}
Let now $u,v$ be arbitrary dyadic numbers such that $m-1\leqslant
u\leqslant v\leqslant m$. Since
\[
 \mathbb{E}_{t_{k-1}^{(n)}}(|X(t_k^{(n)})|^2-|X(t_{k-1}^{(n)})|^2)=
 \mathbb{E}_{t_{k-1}^{(n)}}(|X(t_k^{(n)})-X(t_{k-1}^{(n)})|^2)\geqslant0,
\]
we obtain that for sufficiently large $n$ (again so large that
\linebreak $u,v\in\{t_k^{(n)}\colon k=(m-1)2^n+1,\dots,m2^n\}$)
\[
 S_n^{(m)}(u)\leqslant S_n^{(m)}(v),
\]
and passing to the limit in the above inequality along the
sequence $\{n_r\}$ yields on account of \eqref{e7}
\[
 A(u)-(A_1+\cdots+A_{m-1})\leqslant A(v)-(A_1+\cdots+A_{m-1}),
\]
so
\[
 A(u)\leqslant A(v).
\]
The right continuity of the process $(A(t))$ implies that it is
increasing in $[m-1,m]$, thus from the arbitrariness of
$m,\;(A(t))$ is increasing on the whole of $[0,+\infty)$.
Moreover, $(A(t))$ is positive since $A(0)=0$. The equality
\eqref{e6} gives
\[
 |X(t)|^2=\mathbb{E}_t(|X(m)|^2-(A_1+\cdots+A_m))+A(t)=M(t)+A(t),
\]
for $t\in[m-1,m]$, where
\[
 M(t)=\mathbb{E}_t(|X(m)|^2-(A_1+\cdots+A_m)).
\]
Take an arbitrary $t\in[0,+\infty)$, and let $m$ be a positive
integer such that $t\in[m-1,m]$. For $m-1\leqslant s\leqslant t$
we obviously have $\mathbb{E}_sM(t)=M(s)$. Assume that
$s\in[m-2,m-1]$. Then
\begin{align*}
 \mathbb{E}_sM(t)&=\mathbb{E}_s(|X(m)|^2-(A_1+\cdots+A_m))\\
 &=\mathbb{E}_s|X(m)|^2-\mathbb{E}_sA_m-\mathbb{E}_s(A_1+\cdots+A_{m-1}),\\
\intertext{and}
 M(s)&=\mathbb{E}_s(|X(m-1)|^2-(A_1+\cdots+A_{m-1}))\\
 &=\mathbb{E}_s|X(m-1)|^2-\mathbb{E}_s(A_1+\cdots+A_{m-1}),
\end{align*}
so we obtain
\[
 \mathbb{E}_sM(t)-M(s)=\mathbb{E}_s|X(m)|^2-
 \mathbb{E}_sA_m-\mathbb{E}_s|X(m-1)|^2.
\]
Applying $\mathbb{E}_s$ to both sides of \eqref{e5} yields
\[
 \mathbb{E}_sA_m=\mathbb{E}_s|X(m)|^2-\mathbb{E}_s|X(m-1)|^2,
\]
which shows that
\[
 \mathbb{E}_sM(t)-M(s)=0,\qquad\text{i.e.}\qquad\mathbb{E}_sM(t)=M(s).
\]
Now for an arbitrary $s<t$ choose $s_1<\dots<s_l$ between $s\text{
and }t$ lying in neighbouring intervals with the ends being
positive integers. Then
\begin{align*}
 \mathbb{E}_sM(t)&=\mathbb{E}_s\mathbb{E}_{s_1}\dots\mathbb{E}_{s_l}M(t)
 =\mathbb{E}_s\mathbb{E}_{s_1}\dots\mathbb{E}_{s_{l-1}}M(s_l)\\
 &=\dots=\mathbb{E}_sM(s_1)=M(s),
\end{align*}
proving that $(M(t)\colon t\in[0,+\infty))$ is a martingale.
\end{proof}
Let us now consider an important question concerning the
uniqueness of a Doob-Meyer decomposition. Recall the following
definition from \cite{BW1} (Definition 2.2).
\begin{definition}
An $L^1(\mathcal{A})$ process $(A(t)\colon t\in[0,+\infty))$ is
natural if for each $t>0$ and any sequence $\{\theta_n\}$ of
partitions of $[0,t]$,\\
$\theta_n=\{0=t_0^{(n)}<t_1^{(n)}\dots<t_{m_n}^{(n)}=t\}$ with
$\|\theta_n\|\to 0$ we have
\begin{equation}\label{n}
 \lim_{n\to\infty}\tau\Bigl(\sum_{k=1}^{m_n}
 \mathbb{E}_{t_{k-1}^{(n)}}(y)(A(t_k^{(n)})-A(t_{k-1}^{(n)}))\Bigr)=
 \tau(yA(t)),
\end{equation}
for all $y\in\mathcal{A}$.
\end{definition}
It was pointed out in \cite{BW1} that, in full analogy with the
classical case, if the process $(A(t))$ in decomposition
\eqref{DM} is natural then this decomposition is unique. Note that
equality \eqref{n} may be, on account of \eqref{ce}, rewritten in
the following form (this was also observed in \linebreak
\cite[Definition 5.3]{BW3})
\begin{align*}
 &\lim_{n\to\infty}\tau\Bigl(\sum_{k=1}^{m_n}
 \mathbb{E}_{t_{k-1}^{(n)}}(y)(A(t_k^{(n)})-A(t_{k-1}^{(n)}))\Bigr)\\
 =&\lim_{n\to\infty}\sum_{k=1}^{m_n}
 \tau(\mathbb{E}_{t_{k-1}^{(n)}}(y)(A(t_k^{(n)})-A(t_{k-1}^{(n)})))\\
 =&\lim_{n\to\infty}\sum_{k=1}^{m_n}
 \tau(y\mathbb{E}_{t_{k-1}^{(n)}}((A(t_k^{(n)})-A(t_{k-1}^{(n)}))))\\
 =&\lim_{n\to\infty}\tau \Bigl(y\sum_{k=1}^{m_n}
 \mathbb{E}_{t_{k-1}^{(n)}}(A(t_k^{(n)})-A(t_{k-1}^{(n)}))\Bigr)=\tau(yA(t)),
\end{align*}
which simply means that
\begin{equation}\label{n1}
 \sum_{k=1}^{m_n}\mathbb{E}_{t_{k-1}^{(n)}}(A(t_k^{(n)})-A(t_{k-1}^{(n)}))
 \to A(t) \qquad \text{\em weakly}.
\end{equation}
From decomposition \eqref{DM} it follows that for any $0\leqslant
s\leqslant t$ we have
\begin{align*}
 \mathbb{E}_s(A(t)-A(s))&=\mathbb{E}_s(|X(t)|^2-|X(s)|^2)
 -\mathbb{E}_s(M(t)-M(s))\\&=\mathbb{E}_s(|X(t)|^2-|X(s)|^2),
\end{align*}
by martingale property, so condition \eqref{n1} becomes
\begin{equation}\label{n2}
 \sum_{k=1}^{m_n}\mathbb{E}_{t_{k-1}^{(n)}} (|X(t_k^{(n)})|^2-|X(t_{k-1}^{(n)})|^2)\to
 A(t)\qquad \text{\em weakly}.
\end{equation}
In the next section we shall show that this condition is satisfied
for a certain class of martingales, thus obtaining the uniqueness
of a Doob-Meyer decomposition. Moreover, an explicit form of this
decomposition will be given.

\section{Quadratic variation process}
In this section we assume that $(X(t)\colon t\in[0,+\infty))$ is
an $\mathcal{A}$-valued martingale continuous in Segal's sense.
The following theorem is a noncommutative counterpart of
Theorem~4.1 from \cite{CW}.
\begin{theorem}\label{qv0}
Let $t\geqslant0$, and denote by\\
$\theta=\{0=t_0<t_1<\dots<t_m=t\}$ a partition of $[0,t]$. Put
\[
 S_{\theta}(t)=\sum_{k=1}^m|\Delta X(t_k)|^2.
\]
Then $\{S_{\theta}(t)\}$ converges in $\|\cdot\|_2$-norm as
$\|\theta\|\to 0$ to
\begin{equation}\label{qv}
 \begin{aligned}
 \langle X\rangle(t):=&|X(t)|^2-|X(0)|^2+\\
 -&\left(\int_0^t dX^*(u)\,X(u)
 +\int_0^t X^*(u)\,dX(u)\right).
 \end{aligned}
\end{equation}
\end{theorem}
\begin{proof}
We have, analogously as in the classical case,
\begin{align*}
 S_{\theta}(t)&=\sum_{k=1}^m[X(t_k)-X(t_{k-1})]^*
 [X(t_k)-X(t_{k-1})]\\&=\sum_{k=1}^m
 (|X(t_k)|^2-|X(t_{k-1})|^2)+\\
 &-\sum_{k=1}^m[X^*(t_k)-X^*(t_{k-1})]X(t_{k-1})
 +\\&-\sum_{k=1}^mX^*(t_{k-1})[X(t_k)-X(t_{k-1})]\\
 &=|X(t)|^2-|X(0)|^2-\sum_{k=1}^m
 [X^*(t_k)-X^*(t_{k-1})]X(t_{k-1})
 +\\&-\sum_{k=1}^m X^*(t_{k-1})[X(t_k)-X(t_{k-1})].
\end{align*}
Now observe that the martingales $(X(t))\text{ and }(X^*(t))$
satisfy the assumption of Theorem ~\ref{Int} both as the
integrators and the integrands, thus the two sums in the equation
above, being respectively the left and right integral sum, tend to
the integrals $\int_0^t dX^*(u)\,X(u)\text{ and }\int_0^t
X^*(u)\,dX(u)$, which proves the theorem.
\end{proof}
For any $0\leqslant s\leqslant t$ we clearly have $0\leqslant
S_{\theta}(s)\leqslant S_{\theta}(t)$, so the process $(\langle
X\rangle(t)\colon t\in[0,+\infty))$ defined by equation
\eqref{qv}, being the limit of $S_{\theta}(t)$, is positive and
increasing; obviously $\langle X\rangle(0)=0$. This process is
called the \emph{quadratic variation process} for the martingale
$(X(t)\colon t\in[0,+\infty))$. Denote
\[
 M(t)=|X(0)|^2+\int_0^t dX^*(u)\,X(u)
 +\int_0^t X^*(u)\,dX(u).
\]
By virtue of Theorem ~\ref{Int} we obtain that $(M(t))$ is a
martingale, thus \eqref{qv} gives a Doob-Meyer decomposition
\begin{equation}\label{DM1}
 |X(t)|^2=M(t)+\langle X\rangle(t).
\end{equation}
From Theorem~\ref{Scm} it follows that $(M(t))$ is weakly
continuous in Segal's sense. We shall show that the quadratic
variation process $(\langle X\rangle(t))$ is natural. The
following lemma is a noncommutative version of Lemma 5.10 from
\cite{KS}.
\begin{lemma}\label{L}
Let $t>0$, and denote by $\theta=\{0=t_0<t_1<\dots<t_m=t\}$ a
partition of $[0,t]$. Then
\[
 \lim_{\|\theta\|\to 0}\tau\Bigl(\sum_{k=1}^m
 |\Delta X(t_k)|^4\Bigr)=0.
\]
\end{lemma}
\begin{proof}
The estimations in the proof are similar to those in Theorem
~\ref{Int}. Put
\begin{align}
 M&=\|X(t)\|_{\infty}=\sup_{0\leqslant s\leqslant
 t}\|X(s)\|_{\infty},\label{o11}\\
 K&=\tau(|X(t)|^2-|X(0)|^2),\label{o21}
\end{align}
and let $\alpha_4$ be as in \eqref{PX} with $p=4$. Let
$\varepsilon>0$ be given. From Segal's continuity of the
martingale $(X(u))$ there exist a projection $e$ in $\mathcal{A}$
with
\begin{equation}\label{o31}
 \tau(e^{\bot})<\frac{\varepsilon^2}{64M^4\alpha_4^4\|X(t)\|_4^4},
\end{equation}
and $\delta>0$ such that for any $t',t''\in[0,t]\text{ with }
|t'-t''|<\delta$, we have
\begin{equation}\label{e21}
 \|[X(t')-X(t'')]e\|_{\infty}<\frac{\varepsilon}{4MK}.
\end{equation}
Let $\theta=\{0=t_0<t_1<\dots<t_m=t\}$ be a partition of $[0,t]$
such that $\|\theta\|<\delta$. We have
\begin{align*}
 \tau\Bigl(\sum_{k=1}^m|\Delta X(t_k)|^4\Bigr)
 &=\tau\Bigl(\sum_{k=1}^m e|\Delta X(t_k)|^4\Bigr)+
 \tau\Bigl(\sum_{k=1}^m e^{\bot}|\Delta X(t_k)|^4\Bigr)\\&=I_1+I_2,
\end{align*}
where
\begin{equation*}
 I_1=\tau\Bigl(\sum_{k=1}^m e|\Delta X(t_k)|^4\Bigr),\qquad
 I_2=\tau\Bigl(\sum_{k=1}^m e^{\bot}|\Delta X(t_k)|^4\Bigr).
\end{equation*}
For $I_1$ we have using \eqref{o11}, \eqref{o21} and \eqref{e21}
together with \eqref{e1}
\begin{align*}
 |I_1|&=\sum_{k=1}^m |\tau(e|\Delta X(t_k)|^4)|
 \leqslant\sum_{k=1}^m \|e|\Delta X(t_k)|^2\|_{\infty}\||\Delta
 X(t_k)|^2\|_1\\
 &\leqslant\sum_{k=1}^m \|e\Delta
 X(t_k)\|_{\infty}\||\Delta X(t_k)|\|_{\infty}\||\Delta
 X(t_k)|^2\|_1\\
 &<\sum_{k=1}^m 2M\frac{\varepsilon}{4MK}\||\Delta
 X(t_k)|^2\|_1=\frac{\varepsilon}{2K}\sum_{k=1}^m \||\Delta
 X(t_k)|^2\|_1=\frac{\varepsilon}{2}.
\end{align*}
For $I_2$ we have
\begin{align*}
 |I_2|&=\sum_{k=1}^m |\tau(e^{\bot}|\Delta X(t_k)|^4)|\\
 &=\sum_{k=1}^m |\tau(|\Delta X(t_k)|e^{\bot}|\Delta X(t_k)||\Delta
 X(t_k)|^2)|\\
 &\leqslant\sum_{k=1}^m \||\Delta X(t_k)|e^{\bot}|\Delta X(t_k)|\|_1 \||\Delta
 X(t_k)|^2\|_{\infty}\\&=\sum_{k=1}^m \tau(e^{\bot}|\Delta X(t_k)|^2)
 \||\Delta X(t_k)|^2\|_{\infty}\leqslant4M^2\tau\Bigl(\sum_{k=1}^m
 e^{\bot}|\Delta X(t_k)|^2\Bigr)\\&\leqslant4M^2\Big\|e^{\bot}\sum_{k=1}^m
 |\Delta X(t_k)|^2\Big\|_1.
\end{align*}
From H\"{o}lder's inequality we get
\begin{align*}
 &\Big\|e^{\bot}\sum_{k=1}^m|\Delta
 X(t_k)|^2\Big\|_1\leqslant\Big\|e^{\bot}\Big\|_2\Big\|\sum_{k=1}^m|\Delta
 X(t_k)|^2\Big\|_2\\=&[\tau(e^{\bot})]^{1/2}\Big\|\Bigl(\sum_{k=1}^m|\Delta
 X(t_k)|^2\Bigr)^{1/2}\Big\|_4^2,
\end{align*}
and taking into account inequalities \eqref{PX}, \eqref{o31} we
finally obtain
\begin{align*}
 |I_2|&\leqslant4M^2\Big\|e^{\bot}\sum_{k=1}^m
 |\Delta X(t_k)|^2\Big\|_1\leqslant4M^2[\tau(e^{\bot})]^{1/2}\Big\|\Bigl(\sum_{k=1}^m
 |\Delta X(t_k)|^2\Bigr)^{1/2}\Big\|_4^2\\
 &<4M^2\frac{\varepsilon}{8M^2\alpha_4^2\|X(t)\|_4^2}
 \alpha_4^2\|X(t)\|_4^2=\frac{\varepsilon}{2},
\end{align*}
and thus
\[
 \tau\Bigl(\sum_{k=1}^m|\Delta
 X(t_k)|^4\Bigr)=I_1+I_2<\varepsilon,
\]
which finishes the proof.
\end{proof}
\begin{theorem}
The process $(\langle X\rangle(t)\colon t\in[0,+\infty))$ is
natural.
\end{theorem}
\begin{proof}
Let $t>0$, and denote by $\theta=\{0=t_0<t_1<\dots<t_m=t\}$ a
partition of $[0,t]$. We shall first show that
\[
 \lim_{\|\theta\|\to 0}\Big\|\sum_{k=1}^m(|\Delta
 X(t_k)|^2-\mathbb{E}_{t_{k-1}}|\Delta X(t_k)|^2)\Big\|_2=0.
\]
We have
\begin{align*}
 &\Big\|\sum_{k=1}^m(|\Delta
 X(t_k)|^2-\mathbb{E}_{t_{k-1}}|\Delta X(t_k)|^2)\Big\|_2^2\\
 =&\tau\Bigl(\sum_{i,k=1}^m(|\Delta
 X(t_i)|^2-\mathbb{E}_{t_{i-1}}|\Delta X(t_i)|^2)(|\Delta
 X(t_k)|^2-\mathbb{E}_{t_{k-1}}|\Delta X(t_k)|^2)\Bigr)\\
 =&\sum_{i,k=1}^m\tau((|\Delta
 X(t_i)|^2-\mathbb{E}_{t_{i-1}}|\Delta X(t_i)|^2)(|\Delta
 X(t_k)|^2-\mathbb{E}_{t_{k-1}}|\Delta X(t_k)|^2)).
\end{align*}
For $i<k$ we obtain
\begin{align*}
 &\tau((|\Delta
 X(t_i)|^2-\mathbb{E}_{t_{i-1}}|\Delta X(t_i)|^2)(|\Delta
 X(t_k)|^2-\mathbb{E}_{t_{k-1}}|\Delta X(t_k)|^2))\\
 =&\tau(\mathbb{E}_{t_{k-1}}(|\Delta
 X(t_i)|^2-\mathbb{E}_{t_{i-1}}|\Delta X(t_i)|^2)(|\Delta
 X(t_k)|^2-\mathbb{E}_{t_{k-1}}|\Delta X(t_k)|^2))\\
 =&\tau((|\Delta
 X(t_i)|^2-\mathbb{E}_{t_{i-1}}|\Delta X(t_i)|^2)\mathbb{E}_{t_{k-1}}(|\Delta
 X(t_k)|^2-\mathbb{E}_{t_{k-1}}|\Delta X(t_k)|^2))\\=&0,
\end{align*}
and analogously for $i>k$. Consequently, by Lemma ~\ref{L}
\begin{align*}
 &\Big\|\sum_{k=1}^m(|\Delta
 X(t_k)|^2-\mathbb{E}_{t_{k-1}}|\Delta X(t_k)|^2)\Big\|_2^2\\
 =&\tau\Bigl(\sum_{k=1}^m\bigl[|\Delta
 X(t_k)|^2-\mathbb{E}_{t_{k-1}}|\Delta X(t_k)|^2\bigr]^2\Bigr)\\
 \leqslant&\tau\Bigl(\sum_{k=1}^m 2\bigl[|\Delta
 X(t_k)|^4+\bigl(\mathbb{E}_{t_{k-1}}|\Delta X(t_k)|^2\bigr)^2\bigr]\Bigr)\\
 \leqslant&2\tau\Bigl(\sum_{k=1}^m \bigl[|\Delta
 X(t_k)|^4+\mathbb{E}_{t_{k-1}}|\Delta X(t_k)|^4\bigr]\Bigr)\\
 =&4\tau\Bigl(\sum_{k=1}^m|\Delta X(t_k)|^4\Bigr)\to 0,
\end{align*}
where in the last estimation we have used the inequality
\[
 (\mathbb{E}_tx)^2\leqslant\mathbb{E}_tx^2
\]
valid for $x=x^*$ because $\mathbb{E}_t$ is positive.

From Theorem ~\ref{qv0} we have
\[
 \lim_{\|\theta\|\to0}\sum_{k=1}^m|\Delta X(t_k)|^2=
 \langle X\rangle(t),
\]
which yields
\[
 \lim_{\|\theta\|\to 0}\sum_{k=1}^m\mathbb{E}_{t_{k-1}}|\Delta
 X(t_k)|^2=\langle X\rangle(t)\qquad \text{\em
 in }\|\cdot\|_2-\text{\em norm},
\]
and thus condition \eqref{n2} is satisfied.
\end{proof}
Our next aim is to show a noncommutative counterpart of the
classical result saying that if the martingale in the Doob-Meyer
decomposition is continuous with probability one then this
decomposition is unique.
\begin{theorem}\label{DMu}
Let
\[
 |X(t)|^2=M(t)+A(t)
\]
be a Doob-Meyer decomposition of the submartingale $(|X(t)|^2)$
such that the martingale $(M(t))$ is weakly continuous in Segal's
sense and $M(0)=|X(0)|^2$. Then this decomposition is unique.
\end{theorem}
\begin{proof}
Assume that there are decompositions
\[
 |X(t)|^2=M_1(t)+A_1(t)=M_2(t)+A_2(t),
\]
where $(M_1(t)),\:(M_2(t))$ are weakly continuous in Segal's sense
martingales such that $M_1(0)=M_2(0)=|X(0)|^2,\text{ and
}(A_1(t)),\:(A_2(t))$ are increasing positive processes. Put
\[
 M(t)=M_1(t)-M_2(t).
\]
Then $(M(t))$ is a weakly continuous in Segal's sense martingale
such that
\[
 M(t)=A_2(t)-A_1(t),
\]
in particular, $(M(t))$ is selfadjoint and $M(0)=0$. Fix $t>0$.
Take an arbitrary $\varepsilon>0$. From weak uniform continuity in
Segal's sense of the martingale $(M(s)\colon s\in[0,t])$ it
follows that there exist a projection $e\in\mathcal{A}$ with
\begin{equation}\label{o32}
 \tau(e^{\bot})<\frac{\varepsilon^2}{16\alpha_8^4\|M(t)\|_8^4},
\end{equation}
where $\alpha_8$ is as in \eqref{PX} for $p=8$, and $\delta>0$
such that for any \linebreak $t',t''\in[0,t]\text{ with
}|t'-t''|<\delta$ we have
\begin{equation}\label{e22}
 \|e[M(t')-M(t'')]e\|_{\infty}\leqslant\frac{\varepsilon}
 {4(\|A_1(t)\|_1+\|A_2(t)\|_1)}.
\end{equation}
Let $0=t_0<t_1<\dots<t_m=t$ be a partition of the interval $[0,t]$
such that $\max_{1\leqslant k\leqslant m}(t_k-t_{k-1})<\delta$. We
have
\begin{align*}
 &\tau\Big(\sum_{k=1}^m[\Delta M(t_k)]^2\Big)=\tau\Big(\sum_{k=1}^m
 e\Delta M(t_k)e\Delta M(t_k)\Big)+\\ &\tau\Big(\sum_{k=1}^m
 e^{\bot}\Delta M(t_k)e\Delta M(t_k)\Big)+\tau\Big(\sum_{k=1}^m
 e\Delta M(t_k)e^{\bot}\Delta M(t_k)\Big)+\\ &\tau\Big(\sum_{k=1}^m
 e^{\bot}\Delta M(t_k)e^{\bot}\Delta M(t_k)\Big)=I_1+I_2+I_3+I_4,
\end{align*}
where
\begin{align*}
 I_1&=\tau\Big(\sum_{k=1}^m e\Delta M(t_k)e\Delta M(t_k)\Big),&
 I_2&=\tau\Big(\sum_{k=1}^m e^{\bot}\Delta M(t_k)e\Delta
 M(t_k)\Big),\\ I_3&=\tau\Big(\sum_{k=1}^m e\Delta
 M(t_k)e^{\bot}\Delta M(t_k)\Big),& I_4&=\tau\Big(\sum_{k=1}^m
 e^{\bot}\Delta M(t_k)e^{\bot}\Delta M(t_k)\Big).
\end{align*}
On account of \eqref{e22} we get
\begin{align*}
 |I_1|&=I_1=\sum_{k=1}^m|\tau(e\Delta M(t_k)e\Delta M(t_k)|\leqslant
 \sum_{k=1}^m\|e\Delta M(t_k)e\|_{\infty}\|\Delta M(t_k)\|_1\\
 &\leqslant\frac{\varepsilon}{4(\|A_1(t)\|_1+\|A_2(t)\|_1)}
 \sum_{k=1}^m\|\Delta M(t_k)\|_1,
\end{align*}
and since $(A_1(t))\text{ and }(A_2(t))$ are increasing and
$A_1(0)=A_2(0)=0$,
\begin{align*}
 &\sum_{k=1}^m\|\Delta M(t_k)\|_1= \sum_{k=1}^m\|\Delta A_2(t_k)-
 \Delta A_1(t_k)\|_1\\ \leqslant&\sum_{k=1}^m\|\Delta A_2(t_k)\|_1+
 \sum_{k=1}^m\|\Delta A_1(t_k)\|_1\\
 =&\sum_{k=1}^m\tau(A_2(t_k)-A_2(t_{k-1}))
 +\sum_{k=1}^m\tau(A_1(t_k)-A_1(t_{k-1}))\\
 =&\tau(A_2(t))+\tau(A_1(t))=\|A_2(t)\|_1+\|A_1(t)\|_1,
\end{align*}
which gives the estimation
\[
 I_1\leqslant\frac{\varepsilon}{4}.
\]
Furthermore on account of \eqref{o32} we obtain
\begin{align*}
 |I_2|&=|I_3|=I_2=I_3=\Big|\tau\Big(\sum_{k=1}^me^{\bot}\Delta M(t_k)e\Delta
 M(t_k)\Big)\Big|\\ &=\Big|\tau\Big(e^{\bot}\sum_{k=1}^m\Delta M(t_k)e\Delta
 M(t_k)\Big)\Big|\leqslant\Big\|e^{\bot}\Big\|_2\Big\|\sum_{k=1}^m\Delta
 M(t_k)e\Delta M(t_k)\Big\|_2\\ &=\Big\|e^{\bot}\Big\|_2\Big\|\sum_{k=1}^m(\Delta
 M(t_k)e)(\Delta M(t_k)e)^*\Big\|_2\\ &=[\tau(e^{\bot})]^{1/2}\Big\|\Big(\sum_{k=1}^m
 (\Delta M(t_k)e)(\Delta M(t_k)e)^*\Big)^{1/2}\Big\|_4^2\\
 &<\frac{\varepsilon}{4\alpha_8^2\|M(t)\|_8^2}\Big\|\Big(\sum_{k=1}^m
 (\Delta M(t_k)e)(\Delta M(t_k)e)^*\Big)^{1/2}\Big\|_4^2.
\end{align*}
From \cite[Lemma~1.1]{PX} it follows that
\[
 \Big\|\Big(\sum_{k=1}^m(\Delta M(t_k)e)(\Delta M(t_k)e)^*\Big)^
 {1/2}\Big\|_4\leqslant\Big\|\Big(\sum_{k=1}^m
 [\Delta M(t_k)]^2\Big)^{1/2}\Big\|_8\Big\|e\Big\|_8,
\]
and from inequality \eqref{PX} we obtain
\[
 \Big\|\Big(\sum_{k=1}^m[\Delta M(t_k)]^2\Big)^{1/2}\Big\|_8
 \leqslant\alpha_8\|M(t)\|_8.
\]
Consequently,
\begin{align*}
 I_2=I_3&<\frac{\varepsilon}{4\alpha_8^2\|M(t)\|_8^2}\Big\|\Big(\sum_{k=1}^m
 (\Delta M(t_k)e)(\Delta M(t_k)e)^*\Big)^{1/2}\Big\|_4^2\\ &\leqslant
 \frac{\varepsilon}{4\alpha_8^2\|M(t)\|_8^2}\Big\|\Big(\sum_{k=1}^m
 [\Delta M(t_k)]^2\Big)^{1/2}\Big\|_8^2\Big\|e\Big\|_8^2\\
 &\leqslant\frac{\varepsilon}{4\alpha_8^2\|M(t)\|_8^2}
 \alpha_8^2\|M(t)\|_8^2=\frac{\varepsilon}{4}.
\end{align*}
The same estimation is obtained for $|I_4|=I_4$, so finally we
have
\[
 \tau\Big(\sum_{k=1}^m[\Delta M(t_k)]^2\Big)\leqslant
 I_1+I_2+I_3+I_4<\varepsilon.
\]
On the other side, from equality \eqref{e1} and martingale
property we get
\begin{align*}
 \tau(|M(t)-M(0)|^2)=&\tau(|M(t)|^2-|M(0)|^2)=\sum_{k=1}^m\||\Delta
 M(t_k)|^2\|_1\\ =&\tau\Big(\sum_{k=1}^m[\Delta M(t_k)]^2\Big)
 <\varepsilon,
\end{align*}
and since $\varepsilon>0$ was arbitrary this yields
\[
 \tau(|M(t)-M(0)|^2)=0.
\]
The faithfulness of $\tau$ implies that
\[
 M(t)=M(0)=0,
\]
which gives the uniqueness of the decomposition.
\end{proof}
\begin{remark}
From the decomposition \eqref{DM1} and the uniqueness it follows
that the process $(A(t))$ in the above theorem is actually equal
to $(\langle X\rangle(t))$.
\end{remark}
We finish with a brief discussion of the classical notion of the
cross-variation process.
\begin{definition}
Let $(X(t)\colon t\in[0,+\infty))\text{ and }(Y(t)\colon
t\in[0,+\infty))$ be \linebreak $\mathcal{A}$-valued continuous in
Segal's sense martingales. The \emph{cross-variation process}
(called also the \emph{bracket} or \emph{quadratic covariation
process}) is defined by
\[
\begin{aligned}
 \langle X,Y\rangle(t)&=\frac{1}{4}\{\langle X+Y\rangle(t)-\langle
 X-Y\rangle(t)\\&+i[\langle iX+Y\rangle(t)-\langle iX-Y\rangle(t)]\}.
\end{aligned}
\]
\end{definition}
This is of course the polarization formula for the `sesquilinear
form' (linear in the second position) $\langle X,Y\rangle$
obtained from the `quadratic form' $\langle X\rangle$, so that
\[
 \langle X,X\rangle(t)=\langle X\rangle(t).
\]
After straightforward calculations we obtain
\begin{align*}
 \langle
 X,Y\rangle(t)&=X(t)^*Y(t)-X(0)^*Y(0)+\\&-\left[\int_0^tdX^*(u)\,Y(u)
 +\int_0^tX^*(u)\,dY(u)\right].
\end{align*}
Our last theorem shows that the cross-variation process can be
obtained in the way entirely analogous to the classical case.
\begin{theorem}
Let $(X(t))\text{ and }(Y(t))$ be as above, and let\linebreak
$\theta=\{0=t_0<t_1<\dots<t_m=t\}$ denote a partition of $[0,t]$.
Then
\[
 \langle X,Y\rangle(t)=\lim_{\|\theta\|\to 0}\sum_{k=1}^m
 [X(t_k)-X(t_{k-1})]^*[Y(t_k)-Y(t_{k-1})]
 \text{ in }\|\cdot\|_2-\text{norm}.
\]
\end{theorem}
\begin{proof}
We have
\begin{align*}
 &\sum_{k=1}^m [X(t_k)-X(t_{k-1})]^*[Y(t_k)-Y(t_{k-1})]\\=
 &\sum_{k=1}^m X^*(t_k)Y(t_k)-\sum_{k=1}^m X^*(t_k)Y(t_{k-1})
 -\sum_{k=1}^m X^*(t_{k-1})[Y(t_k)-Y(t_{k-1})]\\=
 &X^*(t)Y(t)-X^*(0)Y(0)+\sum_{k=1}^m
 X^*(t_{k-1})Y(t_{k-1})-\sum_{k=1}^m X^*(t_k)Y(t_{k-1})+\\
 -&\sum_{k=1}^m
 X(t_{k-1})[Y(t_k)-Y(t_{k-1})]=X^*(t)Y(t)-X^*(0)Y(0)+\\
 -&\sum_{k=1}^m [X^*(t_k)-X^*(t_{k-1})]Y(t_{k-1})-\sum_{k=1}^m
 X^*(t_{k-1})[Y(t_k)-Y(t_{k-1})].
\end{align*}
The two sums on the right hand side of the above equation tend to
$\int_0^tdX^*(u)\,Y(u)\text{ and }\int_0^t X^*(u)\,dY(u)$,
respectively, which shows the claim.
\end{proof}

\end{document}